\documentclass{article}
\usepackage{graphicx, amssymb, latexsym, amsfonts, amsmath, lscape, amscd, multirow, multicol,hyperref,
amsthm, color, epsfig, mathrsfs, tikz, enumerate, soul, subcaption, scrextend, float}
\usepackage[colorinlistoftodos,backgroundcolor=orange!40,linecolor=orange]{todonotes}

\usetikzlibrary{calc}

\theoremstyle{definition}

\newtheorem{theorem}{Theorem}
\newtheorem{conjecture}[theorem]{Conjecture}

\newtheorem{proposition}[theorem]{Proposition}

\newtheorem{observation}[theorem]{Observation}

\newtheorem{claim}{Claim}
\renewcommand*{\theclaim}{\Alph{claim}}

\newtheorem{case}{Case}

\newtheorem{subcase}{Subcase}[case]

\newtheorem{subclaim}{Subclaim}[claim]

\newcommand{\gLT}{\gamma_t^L}

\newcommand{\smallqed}{{\small ($\Box$)}}
\newcommand{\smallqedbis}{{\small ($\diamond$)}}
\newcommand{\tdom}{{\rm tdom}}
\newcommand{\cF}{{\cal F}}

\newenvironment{claimproof}{\noindent\emph{Proof of claim.}}{~\smallqed\newline\medskip}
\newenvironment{subclaimproof}{\noindent\emph{Proof of subclaim.}}{~\smallqedbis\newline\medskip}

\makeatletter

\renewcommand\@makefntext[1]{%
\setlength\parindent{1em}%
\noindent
\makebox[1.5em][r]{\@thefnmark.~}{#1}}

\makeatother

\title{Progress towards the two-thirds conjecture on \\ locating-total dominating sets\footnote{The first two authors were supported by the French government IDEX-ISITE initiative CAP 20-25 (ANR-16-IDEX-0001), the International Research Center ``Innovation Transportation and Production Systems'' of the I-SITE CAP 20-25, and the ANR project GRALMECO (ANR-21-CE48-0004). The research of Anni Hakanen was supported by Jenny and Antti Wihuri Foundation and Academy of Finland, grant number 338797. The research of Michael Henning was supported in part by the South African National Research Foundation, Grant Numbers 132588  and 129265, and the University of Johannesburg.}}

\author{Dipayan Chakraborty\footnote{\noindent Universit\'e Clermont Auvergne, CNRS, Mines de Saint-\'Etienne, Clermont Auvergne INP, LIMOS, 63000 Clermont-Ferrand, France} \footnotemark[4] \and Florent Foucaud\footnotemark[2] \and Anni Hakanen\footnotemark[2] \footnote{Department of Mathematics and Statistics, University of Turku, Finland} \and Michael A. Henning\footnote{\noindent Department of Mathematics and Applied Mathematics, University of Johannesburg} \and Annegret K. Wagler\footnotemark[2]}

\begin{document}

\maketitle
\noindent

\bigskip

\begin{abstract}
We study upper bounds on the size of optimum locating-total dominating sets in graphs. A set $S$ of vertices of a graph $G$ is a locating-total dominating set if every vertex of $G$ has a neighbor in $S$, and if any two vertices outside $S$ have distinct neighborhoods within $S$. The smallest size of such a set is denoted by $\gLT(G)$. It has been conjectured that $\gLT(G)\leq\frac{2n}{3}$ holds for every twin-free graph $G$ of order $n$ without isolated vertices. We prove that the conjecture holds for cobipartite graphs, split graphs, block graphs and subcubic graphs.
\end{abstract}

\section{Introduction}

Our aim is to study upper bounds on the smallest
size of locating-total dominating sets in graphs. This notion is part of the extended research area of \emph{identification problems} in graphs and, more generally, discrete structures like hypergraphs. In these types of problems, one seeks to select a small solution set, generally some vertices of a graph, in order to uniquely identify each vertex of the graph by its relationship with the selected vertices. More precisely, given a set $D$ of vertices of a graph $G$, we say that two vertices $v$ and $w$ of $G$ are \emph{located} by $D$ if they have distinct sets of neighbors in $D$. On the other hand, any set $D$ of vertices of $G$ that locates every pair $v,w$ of vertices in $V(G) \setminus D$ is called a \emph{locating set} of $G$. Various notions based on this location property have been studied, such as locating-dominating sets~\cite{S88}, identifying codes~\cite{KCL98} or separating sets~\cite{BS07}, to name a few. We refer to the online bibliography on these topics maintained by Jean and Lobstein~\cite{biblio} (almost 500 references by the time of writing, in 2022). Such problems have a wide range of applications, such as fault-detection in sensor or computer networks~\cite{KCL98,UTS04}, biological testing~\cite{MS85}, machine learning~\cite{BGL19}, or canonical representations of graphs~\cite{B80,KPSV04}, to name a few.

In this paper, we study the notion of a locating-total dominating set. A set $D$ of vertices is a \emph{total dominating set}, abbreviated TD-set, of a graph $G$ if every vertex in $G$ has a neighbor in $D$. Total dominating sets are a natural and widely studied variant of the domination problem in graphs. We refer to the books~\cite{HaHeHe-23,bookTD} for an overview on the topic. A \emph{locating-total dominating set} of $G$ is a TD-set $D \subset V(G)$ such that any two vertices of $G$ not in $D$ are located by $D$. The smallest size of such a locating-total dominating set of a graph $G$ is called the \emph{locating-total domination number} of $G$ and is denoted by $\gLT(G)$. A graph admits a (locating-)total dominating set if and only if it has no isolated vertex. We abbreviate a locating-total dominating set by LTD-set, and we say that a $\gamma_t^L$-\emph{set of $G$} is an LTD-set of minimum cardinality, $\gamma_t^L(G)$, in $G$.

The concept of a locating-total dominating set was first considered in~\cite{hhh06}, based on the similar concept of a locating-dominating set (where total domination is replaced with usual domination) introduced by Slater in the 1980s~\cite{S88}. It was studied for example in~\cite{ABLW19,ABLW20,BCMMS07,BD11,BFL08,C08,CR09,CS11,conjpaperLTD,line,hl12,hr12}. The associated decision problem is NP-hard~\cite{MRJRM17,PF23,RRS19}. The related concept where \emph{all} vertices (not just the ones outside of $D$) must be located was studied in~\cite{C08,FL22,hhh06}.

It is known that any graph of order~$n$ with every component of order at least~$3$ has a TD-set size at most $\frac{2}{3}n$~\cite{CDH80}, and this bound is tight only for the triangle, the 6-cycle and the family of 2-coronas of graphs~\cite{BCV00}. (The \emph{$2$-corona} $H \circ P_2$ of a connected graph $H$ is the graph of order~$3|V(H)|$ obtained from $H$ by attaching a path of length~$2$ to each vertex of $H$ so that the resulting paths are vertex-disjoint.)

However, such a bound does not hold for locating-total dominating sets. Two vertices of a graph are \emph{twins} if they either have the same open neighborhood (\emph{open twins}) or the same closed neighborhood (\emph{closed twins}). Consider a set $S$ of vertices that are pairwise twins of size at least~2 in a graph $G$ (and $S$ forms either a clique or an independent set). 
Then, any locating-total dominating set $D$ needs to contain all vertices of $S$ except possibly one. Indeed, any two such vertices not in $D$ would otherwise not be located. For example, any complete graph of order at least~2 has only twins, and thus has its locating-total domination number equal to its order minus one (while any two vertices form a total dominating set). Other families of (twin-free) graphs with a total dominating set of size~2 and arbitrarily large locating-total domination number have been described in~\cite{conjpaperLTD}.

Nevertheless, it seems that in the absence of twins, the locating-total domination number cannot be as close to the graph's order as in the general case. Towards such a fact, inspired by a similar problem for (non-total) locating-dominating sets~\cite{cubic,Heia,conjpaper}, two of the authors posed the following conjecture (a graph is called \emph{twin-free} if it does not contain any twins). A graph $G$ is \emph{isolate}-\emph{free} if it contains no isolated vertices.

\begin{conjecture}[\cite{conjpaperLTD}]\label{conj}
Every twin-free isolate-free graph $G$ of order $n$ satisfies $\gLT(G)\leq\frac{2}{3}n$.
\end{conjecture}

Conjecture~\ref{conj} was proved in~\cite{conjpaperLTD} for graphs with no $4$-cycles as subgraphs. It was also proved for line graphs in~\cite{line}. It was also proved in~\cite{conjpaperLTD} to hold for all graphs with minimum degree at least~26 for which another related conjecture~\cite{Heia,conjpaper} holds (which is the case for example for bipartite graphs and cubic graphs). It was proved in a stronger form for claw-free cubic graphs in~\cite{hl12} (there, the $\frac{2}{3}$ factor in the upper bound is in fact replaced with $\frac{1}{2}$, and the authors conjectured that $\frac{1}{2}$ holds for all connected cubic graphs, except $K_4$ and $K_{3,3}$). An approximation of the conjecture was proved to hold for all twin-free graphs in~\cite{conjpaperLTD}, where the $\frac{2}{3}$ factor in the upper bound is replaced with $\frac{3}{4}$.

Note that, if true, the bound of Conjecture~\ref{conj} is tight for the 6-cycle and $2$-coronas, by the following.

\begin{observation}\label{obs:conj-tight}
If a graph $G$ of order $n$ is a triangle, a 6-cycle or a $2$-corona of any graph, then $\gLT(G)=\frac{2}{3}n$.
\end{observation}

A graph is \emph{cobipartite} if its vertex set can be partitioned into two cliques, and \emph{split} if it can be partitioned into a stable set (also called an independent set in the literature) and a clique. A graph is a \emph{block graph} if every $2$-connected component forms a clique. A graph is \emph{subcubic} if each vertex has degree at most~3.

In this paper, we give further evidence towards Conjecture~\ref{conj}, by showing that it holds 
for cobipartite graphs (Section~\ref{sec:cobip}), split graphs (Section~\ref{sec:split}), block graphs (Section~\ref{sec:block}) and subcubic graphs (Section~\ref{sec:subcubic}).
We conclude in Section~\ref{sec:conclu}.

Some of our results are actually slightly stronger. Indeed, the proved upper bound for cobipartite graphs is in fact $\frac{n}{2}$ (which is tight). For twin-free split graphs, we show that the $\frac{2n}{3}$ bound of the conjecture can never be reached. However, we construct infinitely many connected split graphs that come very close to the bound; this is interesting in its own right, showing that not only $2$-coronas have such large locating-total domination numbers. Moreover, the bound for subcubic graphs is proved to hold for all subcubic graphs (except for some small ones like $K_1$, $K_2$, $K_4$ and $K_{1,3}$), even if they have twins.

We now introduce some of the notations used in the paper. The \emph{open} and \emph{closed neighborhoods} of a vertex $v$ in a graph $G$ are denoted $N_G(v)$ and $N_G[v]$, respectively (or $N(v)$ and $N[v]$ if $G$ is clear from the context). We denote by $\deg_G(v)$ the \emph{degree} of the vertex $v$ in the graph $G$, that is, $\deg_G(v) = |N_G(v)|$. The \emph{distance} in $G$ between two vertices $u,v$ is denoted $d_G(u,v)$.
If two graphs $G$ and $H$ are isomorphic, we note $G\cong H$. For a graph $G$ with a vertex or edge $x$, we denote by $G-x$ the subgraph of $G$ obtained by removing $x$, and by $G+x$ the supergraph of $G$ obtained by adding $x$. Similarly, if $X$ is a set of vertices and edges, we use the notations $G-X$ and $G+X$ for the subgraph and supergraph of $G$ obtained by deleting or adding all the elements of $X$. A \emph{leaf} is a vertex of degree~$1$, and its unique neighbor is called a \emph{support vertex}. We denote by $\delta(G)$ and $\Delta(G)$ the minimum and maximum degree, respectively, in the graph $G$. We denote a \emph{path}, a \emph{cycle}, and a \emph{complete} graph on $n$ vertices by $P_n$, $C_n$, and $K_n$, respectively. For an integer $k \ge 1$, we let $[k] = \{1,\ldots,k\}$ and $[k]_0 = \{0,1,\ldots,k\}$.

\section{Cobipartite graphs}\label{sec:cobip}

We now prove a stronger variant of the bound of Conjecture~\ref{conj}, whose proof is a refinement of a similar proof for the (non-total) locating-domination number from~\cite{Heia}. For our next result, we make use of the following result due to Bondy~\cite{B72}.

\begin{theorem}[Bondy~\cite{B72}]\label{lem_bondy}
Let $X$ be a set with $|X| = k$ and let $\mathcal{S} = \{X_1, X_2, \ldots , X_k\}$ be a collection of $k$ distinct subsets of $X$. Then, there exists an element $x$ of $X$ such that $X_i \setminus \{x\} \neq X_j \setminus \{x\}$ for any two sets $X_i, X_j \in \mathcal{S}$ and $i \neq j$. 
\end{theorem}

\begin{theorem}\label{thm:cobip}
For any twin-free cobipartite graph $G$ of order $n$,
we have $\gLT(G) \leq \frac{n}{2}$.
\end{theorem}
\begin{proof}
Let $G$ be a twin-free cobipartite graph of order $n$. If $G$ is disconnected, then $G$ is the disjoint union of two cliques, and thus is not twin-free: it either has closed twins if one of the cliques has order at least~$2$, or is a pair of open twins if both cliques have order~$1$, a contradiction. Hence, the graph $G$ is connected. Let $C_1$ and $C_2$ be two cliques of $G$ that partition its vertex set. Since $G$ is twin-free, both $C_1$ and $C_2$ have size at least~$2$ where we may assume, renaming $C_1$ and $C_2$ if necessary, that $|C_1| \le |C_2|$. Moreover, no two vertices of $C_1$ have the same neighborhood in $C_2$, and vice-versa. Furthermore, at least $|C_1| - 1$ vertices of $C_1$ must have neighbors in $C_2$, and vice-versa. This implies that both $C_1$ and $C_2$ are locating sets of $G$. Thus, if any of $C_1,C_2$ is a TD-set, then it is also an LTD-set. Hence, if $C_1$ is a TD-set, we are done, as $|C_1|\leq\lfloor\frac{n}{2}\rfloor$. If however $C_1$ is not a TD-set, it means that some vertex $v$ of $C_2$ has no neighbors in $C_1$; this vertex is unique since $G$ is twin-free.

If moreover, there is a vertex $w$ in $C_1$ with no neighbor in $C_2$, then we select the set $C=(C_1\setminus\{w\})\cup\{x\}$ as a solution set, where $x\neq v$ is any vertex of $C_2$ other than $v$. This set is clearly a TD-set of $G$. Moreover, any two vertices of $C_2$ are located by $C$, as $N_G(v) \cap C = \{x\}$ and any two other vertices of $C_2\setminus\{x\}$ have distinct and nonempty neighborhoods in $C_1$ (and thus, in $C_1\setminus\{w\}$). Furthermore, $w$ is the only vertex in $V(G) \setminus C$ not dominated by $x$. Hence, $C$ is an LTD-set of $G$. Since $|C|=|C_1|\leq\lfloor\frac{n}{2}\rfloor$, we are done. Therefore, from now on, we assume that every vertex of $C_1$ has a neighbor in $C_2$, more precisely, in the set $C_v = C_2 \setminus \{v\}$. Similarly, if $|C_2| > \lceil \frac{n}{2} \rceil$, that is, if $|C_1|<\lfloor\frac{n}{2}\rfloor$, then by the same preceding arguments, the set $C_1$ together with any vertex of $C_2$ other than $v$ produces an LTD-set of size $|C_1|+1 \leq \lfloor\frac{n}{2}\rfloor$, and we are again done. Hence, we also assume from now on that $|C_2| = \lceil \frac{n}{2} \rceil$.

 Now, we must have $|C_2| \geq 3$, or else, we would have $|C_1| = |C_2| =2$ and by our assumption that every vertex in $C_1$ has a neighbor in $C_2$ and the fact that the vertex $v \in C_2$ has no neighbor in $C_1$, the two vertices of $C_1$ must be twins in $G$, a contradiction. This implies that both $C_2$ and $C_v$ are TD-sets of $G$. Therefore, $C_2$ is an LTD-set of $G$ (recall that $C_2$ is already a locating set of $G$). Therefore, if $n$ is even, then  we have $\gamma^L_t(G) \leq |C_2| = \frac{n}{2}$ and we are done. Hence, for the rest of the proof, we assume that $n$ is odd. Moreover, if $C_v \setminus N_G(z) \neq \emptyset$ for all $z \in C_1$, then we claim that $C_v$ is an LTD-set of $G$. To prove so, since $C_v$ is a TD-set of $G$, we only need to prove that $C_v$ is a locating set of $G$. To begin with, all pairs of vertices of $C_1$ have distinct neighborhoods in $C_v$ and hence, are located by $C_v$. Moreover, with $N_G(v) \cap C_v = C_v$ and the assumption that $C_v \setminus N_G(z) \neq \emptyset$ for all $z \in C_1$, the vertex $v$ is located from every vertex $z$ of $C_1$. This proves the claim that $C_v$ is an LTD-set of $G$. We can therefore assume for the rest of the proof that there exists a vertex $z$ of $C_1$ with $C_v \subset N_G(z)$. Note that there can be at most one such $z \in C_1$ on account of $G$ being twin-free. Next, we prove that there exists a vertex $x \in C_v$ such that the set $C = C_{vx} \cup \{z\}$ is an LTD-set of $G$, where $C_{vx} = C_v \setminus \{x\}$.

Recall that any two vertices in $C_1$ have distinct neighborhoods in $C_v$. Moreover, since $n$ is odd, we have $|C_v| = |C_1| = \lfloor \frac{n}{2} \rfloor$. Therefore, in Theorem~\ref{lem_bondy}, taking $X=C_v$, $k= \lfloor \frac{n}{2} \rfloor$ and $\mathcal{S} = \{N_G(u) \cap C_v : u \in C_1 \}$ as $k$ pairwise distinct sets of $C_v$, by Theorem~\ref{lem_bondy}, there exists a subset $C_{vx}$ of $C_v$, for some $x \in C_v$, which locates all pairs of vertices of $C_1$. In particular, every pair of vertices of $C_1 \setminus \{z\}$ is located by $C_{vx}$. Moreover, since $N_G(z) \cap C_{vx} = C_{vx}$ (using the assumption that $N_G(z) \cap C_v = C_v$), we therefore have $C_{vx} \not \subset N_G(u)$ for all $u \in C_1 \setminus \{z\}$. This implies that the set $C_{vx}$ locates both the vertices $v$ and $x$ from every vertex $u$ of $C_1 \setminus \{z\}$, since $C_{vx} \subset N_G(v)$ and $C_{vx} \subset N_G(x)$. Finally, since the vertex $z$ is a neighbor of $x$ and not of $v$, the vertices $v$ and $x$ are located by $z \in C$. This proves that $C$ is a locating set of $G$. We now show that $C$ is also a TD-set of $G$. To prove so, we see that each vertex of  $C_v \cup C_1 \setminus \{z\}$ has $z \in C$ as its neighbor. Moreover, since $|C_2| \geq 3$, that is, $|C_v| \geq 2$, the vertices $v$ and $z$ have at least one neighbor each in $C_{vx} \subset C$. This proves that $C$ is a TD-set of $G$. Hence, $C$ is an LTD-set of $G$.

Therefore, $|C| = |C_v| = \lfloor \frac{n}{2} \rfloor < \frac{n}{2}$ and we are done again. This proves the theorem.
\end{proof}

The bound of Theorem~\ref{thm:cobip} is tight for \emph{complements of half-graphs} (which are graphs with vertex set $\{x_1,\ldots,x_{2k}\}$ and edge set $\{x_ix_j, |i - j| \leq k-1\}$, see~\cite[Definition~5]{Heia}). These graphs are cobipartite and have their locating(-total) domination number equal to $\frac{n}{2}$~\cite[Proposition~6]{Heia}. More complicated examples can be found in~\cite{Heia}.

\section{Split graphs}\label{sec:split}

Consider a split graph $G=(Q \cup S,E)$ where $Q$ induces a clique and $S$ a stable set.
We suppose that $G$ is isolate-free to ensure the existence of an LTD-set in $G$, which further implies that $G$ is connected and $Q$ non-empty (as every component not containing the clique $Q$ needs to be an isolated vertex from $S$).

\begin{theorem}\label{theorem:split}
  For any twin-free isolate-free split graph $G=(Q \cup S,E)$ of order $n$, we have $\gamma^L_t(G) < \frac{2}{3}n$.
\end{theorem}
\begin{proof}
  First, note that we have $|Q|,|S|\geq 2$ as otherwise $G$ is a single vertex or not twin-free. Therefore, $n \geq 4$.
  Observe next that $Q$ is an LTD-set of $G$ since $Q$ is a TD-set and no two vertices in $S$ have the same neighbors in $Q$ (as $G$ is twin-free) showing that $Q$ is also locating. Hence, the assertion is true if $|Q|< \frac{2}{3}n$, that is, if $|S| > \frac{1}{3}n$. Therefore, we can assume henceforth that $|S| \leq \frac{1}{3}n$. In particular, we can assume that $n \geq 6$ because, otherwise, if $n=5$, we would have $|S| \geq 2 > \frac{1}{3}n$.

  Consider now any set $D$ consisting of all vertices in $S$ and, for each $s \in S$, some arbitrary neighbor $q_s \in Q$ (which exists since $G$ is connected). The set $D$ is an LTD-set of $G$ since $D$ is a TD-set and no two vertices in $Q \setminus D$ have the same neighbors in $S$ (as $G$ is twin-free) implying that $D$ is also locating. Now, we will see how to build such a set $D$ that is also of the required size. Note that there exist two vertices $s,s' \in S$ for which $N(s) \cap N(s') \neq \emptyset$. This is because, if, on the contrary, for each pair of vertices $x,y \in S$, their neighborhoods are disjoint, it implies that
  \begin{enumerate}
  
  \item either the vertices of $N(x)$ are pairwise twins whenever $|N(x)| \geq 2$, a contradiction, or
    \item each set $N(x)$ has cardinality exactly one and so, the rest of the $\frac{1}{3}n$ vertices in $Q$ have no neighbors in $S$. Since $n \geq 6$, it implies that there exist twins in $Q$, again a contradiction. 
  \end{enumerate}
  Thus, let $q_{s,s'} \in N(s) \cap N(s')$ be a common neighbor of $s$ and $s'$. This implies that we can assume the vertices $q_s$ and $q_{s'}$ to be equal to $q_{s,s'}$. This further implies that
  $$|D| = |S| + |\{q_x \in Q : x \in S\} - \{q_{s'}\}| \leq 2|S| - 1 < \frac{2}{3}n.$$
  This proves the result.
\end{proof}

We next show that the bound of Theorem~\ref{theorem:split} cannot be improved for split graphs of orders that are multiples of $3$.

\begin{proposition}\label{prop:split-tight}
For each integer $k\geq 3$, there is a connected twin-free split graph $G_k$ of order $n=3k$ and $\gLT(G_k)= 2k-1$.
\end{proposition}
\begin{proof}
Let $Q = \{q_1, \ldots, q_k\} \cup \{q_1', \ldots, q_k'\}$ be a clique and $S = \{s_1, \ldots, s_k\}$ a stable set, so that 
$N(s_i) = \{q_i,q_i'\}$ for $1 \leq i < k$ and $N(s_k) = \{q_1, \ldots, q_k\}$. 
Note that $q_k'$ has no neighbor in $S$ and that the sets $N(s_i)$ are disjoint for $1 \leq i < k$. See Figure~\ref{fig:split-construction} for an illustration.

\begin{figure}[!htpb]
\centering
\begin{tikzpicture}
\tikzstyle{small node}=[circle, draw, inner sep=0pt, minimum width=6pt]

\node[small node,fill,label={180:$q_1$}](q1) at (0,0)    {};
\node[small node,fill,label={180:$q_1'$}](q1') at (1,0)    {};
\node[small node,label={180:$s_1$}](s1) at (0.5,-2)    {};

\node[small node,fill,label={180:$q_2$}](q2) at (2,0)    {};
\node[small node,fill,label={180:$q_2'$}](q2') at (3,0)    {};
\node[small node,label={180:$s_2$}](s2) at (2.5,-2)    {};

\node[small node,fill,label={180:$q_k$}](qk) at (6,0)    {};
\node[small node,label={180:$q_k'$}](qk') at (7,0)    {};
\node[small node,label={0:$s_k$}](sk) at (6.5,-2)    {};

\node (dots1) at (4.5,0)    {$\ldots$};
\node (dots1) at (4.5,-2)    {$\ldots$};

\node[draw,dashed,rectangle,rounded corners,minimum width=8.5cm,minimum height=1cm,line width=0.5mm] (rQ) at (3.25,0) {};
\node[draw,dashed,rectangle,rounded corners,minimum width=8.5cm,minimum height=1cm,line width=0.5mm] (rS) at (3.25,-2) {};
\node (A) at (8,0)    {$Q$};
\node (A) at (8,-2)    {$S$};

\path[draw,thick]
    (q1) -- (s1) -- (q1')
    (q2) -- (s2) -- (q2')
    (q1) -- (sk) -- (q2) (qk) -- (sk)
;

\end{tikzpicture}\centering
\caption{The construction of graph $G_k$ in the proof of Proposition~\ref{prop:split-tight}, with an optimal LTD-set (black vertices).}\label{fig:split-construction}
\end{figure}
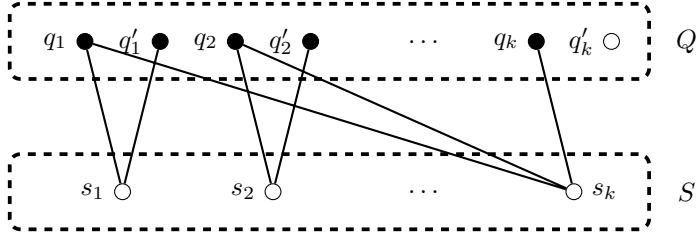

Let $C$ be an LTD-set of $G_k$. 
Consider the $k-1$ closed neighborhoods $N[s_i]$ for $1 \leq i < k$. If we have $|N[s_i]\cap C|\geq 2$ for all $i$ with $1 \leq i < k$, then $\left|\bigcup_{1 \leq i < k} (N[s_i] \cap C) \right| \geq 2k-2$,
and at least one of the remaining vertices $s_k, q_k, q_k'$ must belong to $C$, as otherwise $N(q_k) \cap C = N(q_k') \cap C$ would follow, a contradiction. This implies $|C|\geq 2k-1$.

If, however, for some $i$ with $1 \leq i < k$, we have $|N[s_i]\cap C|=1$, then $s_i\notin C$ since otherwise $s_i$ is not totally dominated by the set $C$. If $N[s_i]\cap C=\{q_i\}$, then $N(q_i')\cap C=Q\cap C$. If $N[s_i]\cap C=\{q_i'\}$, then $N(q_i)\cap C=(Q\cup \{s_k\})\cap C$. 
The two possibilities can occur at most once each. Assume that they both occur once each, with $N[s_a]\cap C=\{q_a'\}$ and $N[s_b]\cap C=\{q_b\}$ (with $1\leq a<b<k$). Note that $s_k\in C$, otherwise $q_a$ and $q_b'$ are not located. Moreover, $C$ must contain $q_k'$ (otherwise $q_b'$ and $q_k'$ are not located) and $q_k$ (otherwise $q_a$ and $q_k$ are not located), and so $|C|\geq 2k-1$, as claimed.

Similarly, if we have $|N[s_i]\cap C|\geq 2$ for all $i$ with $1 \leq i < k$ except that $N[s_a]\cap C=\{q_a'\}$, if $s_k\in C$, then $q_k\in C$, otherwise $q_a$ and $q_k$ are not located. If $s_k\notin C$, then both $q_k,q_k'$ are in $C$ to locate the vertices $q_a,q_k,q_k'$. Thus, again $|C|\geq 2k-1$.

Finally, if we have $|N[s_i]\cap C|\geq 2$ for all $i$ with $1 \leq i < k$ except that $N[s_b]\cap C=\{q_b\}$, if $s_k\in C$, then $q_k'\in C$, otherwise $q_b'$ and $q_k'$ are not located. If $s_k\notin C$, then both $q_k,q_k'$ are in $C$ to locate the vertices $q_b',q_k,q_k'$, and again $|C|\geq 2k-1$.

Thus, in all the above cases, we have $|C|\geq 2k-1$ and, together with the upper bound $\gamma^L_t(G_k) < \frac{2}{3}\times 3k=2k$ from Theorem~\ref{theorem:split}, we finally obtain $\gLT(G_k) = 2k-1$.
\end{proof}


\section{Block graphs}\label{sec:block}

A block graph is a graph in which every maximal $2$-connected subgraph (henceforth referred to as a \emph{block}) is complete. Equivalently, block graphs are diamond-free chordal graphs~\cite{BM86}, where a \emph{diamond} is the graph $K_4-e$, with $e$ being an arbitrary edge of the $K_4$. A \emph{cut-vertex} $v$ of a graph $G$ is one such that the graph $G-v$ has more components than $G$. For any block graph $G$, a \emph{leaf block} of $G$ is a block that contains only a single cut-vertex of $G$. In this section, we show that Conjecture~\ref{conj} holds for block graphs. Trees are a subclass of block graphs in which every block is of order $2$. There are some concepts of trees which we use quite often in proving our result for block graphs and we, therefore, define these concepts formally here.

A \emph{root} of a tree is a fixed vertex of the tree to which the name is designated. Having fixed a root $r$ of a tree $T$, for any vertex $u$ of $T$,

\begin{enumerate}[(1)]
\item a \emph{child} of $u$ is a vertex $v$ of $T$ such that $uv$ is an edge of $T$ and $d_T(v,r)=d_T(u,r)+1$;
\item a \emph{grandchild} of $u$ is a vertex $w$ of $T$ such that $uv$ and $vw$ are edges of $T$ and $d_T(w,r)=d_T(u,r)+2$; and
\item a \emph{great-grandchild} of $u$ is a vertex $x$ of $T$ such that $uv$, $vw$ and $wx$ are edges of $T$ and $d_T(x,r)=d_T(u,r)+3$.
\end{enumerate}

Conversely, the vertex $u$ of $T$ is called the \emph{parent}, the \emph{grandparent} and the \emph{great-grandparent} of $v$, $w$ and $x$, respectively. Given any two vertices $u$ and $v$ of a tree $T$, we say that \emph{$u$ is above $v$ in $T$} (or \emph{$v$ is below $u$ in $T$}) if there exists a sequence of vertices $x_1, x_2, \ldots , x_m$ of $T$ such that, for each $i=1, 2, \ldots , m-1$, $x_{i+1}$ is a child of $x_i$ in $T$, where $m \geq 2$, $x_1 = u$ and $x_m = v$. 

\begin{theorem}\label{theorem_block}
If $G \cong P_3$ or if $G$ is a twin-free isolate-free block graph of order~$n \ge 4$, then $\gamma^L_t(G) \le \frac{2}{3}n$.
\end{theorem}

\begin{proof}
Since the locating-total domination number of a graph is the sum of the locating-total domination numbers of each of the components of the graph, it is therefore enough to prove the theorem for a connected twin-free block graph. Thus, let us assume that $G$ is either isomorphic to a $3$-path or is a connected twin-free block graph of order $n \geq 4$. The proof is by induction on $n \ge 3$. The base case of the induction hypothesis is when $n=3$, in which case $G \cong P_3$ and $\gamma^L_t(G) = 2 = \frac{2}{3}n$. Clearly, any two consecutive vertices of $P_3$ constitute a minimum LTD-set of $P_3$ and hence, the result holds for the base case of the induction hypothesis. We now assume, therefore, that $n \geq 4$ and that the induction hypothesis is true for all connected twin-free
block graphs of order at least~$3$ and at most $n-1$. 
Next, we construct a new graph $T_G$ from $G$ in the following way (see Figure~\ref{block figure} for an example of the construction).

For every block $B$ of $G$, introduce a vertex $u_B \in V(T_G)$ and for every cut-vertex $c \in V(G)$, introduce a vertex $v_c \in V(T_G)$. Next, we introduce edges $u_B v_c \in E(T_G)$ if and only if the cut-vertex $c$ belongs to the block $B$ of $G$.
By construction, therefore, $T_G$ is a tree. Thus, the vertices of the tree $T_G$ are of two types:

\begin{enumerate}[(1)]
\item \emph{$u$-type}: $u_B$ introduced in a one-to-one association with a block $B$ of $G$; and
\item \emph{$v$-type}: $v_c$ introduced in a one-to-one association with a cut-vertex $c$ of $G$.
\end{enumerate}

Notice that any pair of vertices $w,z$ of the tree $T_G$ such that $w$ is the grandparent/grandchild of $z$ in $T_G$ are of the same vertex type. For a fixed cut-vertex $r \in V(G)$, designate $v_r \in V(T_G)$ as the root of $T_G$ (indeed, such a cut-vertex exists as $n \geq 4$ and the twin-free property of $G$ implies that $G$ has at least two blocks). Notice that any leaf of the tree $T_G$ is a vertex of the type $u_B$ for some leaf block $B$ of $G$. By the twin-free nature of $G$, every leaf block $B$ of $G$ has order exactly $2$. Now, fix a leaf $u_{F}$ of $T_G$ that is at the farthest distance, in $T_G$, from the root $v_r$ of $T_G$. We now look at the great-grandparent of the leaf $u_F$ in the tree $T_G$ (indeed, the great-grandparent of $u_F$ in $T_G$ exists because, on account of $G$ being twin-free, at least one of the blocks of $G$ containing the vertex $r$ has a cut-vertex other than $r$).  Notice that the great-grandparent of $u_F$ in $T_G$ must be a vertex of the type $v_p$ for some unique cut-vertex $p$ of $G$. We next define the following.
\begin{flalign*}
&\mathcal{B}_p = \{ B : \text{$B$ is a block of $G$ and $u_B$ is either a child or a great-grandchild of $v_p$ in $T_G$} \};\\
&U = \cup_{B \in \mathcal{B}_p} V(B); \text{ and}\\
&A=\{ x \in U \colon \text{$x$ is a cut-vertex of $G$} \}.
\end{flalign*}


\begin{figure}[t!]
\centering

\begin{subfigure}[h]{0.4\textwidth}
\centering

\begin{tikzpicture}[scale=1.5,
blacknode/.style={circle, draw=black!, fill=black, thick, minimum size= 8mm},
rednode/.style={circle, draw=black!, fill=red!80, thick, minimum size= 8mm},
whitenode/.style={circle, draw=black!, fill=white!, thick, minimum size= 8mm},
]

\node[whitenode, scale=0.4] (a) at (0.3,-1.1) {}; \node at (-0.2,-1) {$r=1$};
\node[whitenode, scale=0.4] (b) at (0.5,-2.1) {}; \node at (0.5,-2.4) {$3$};
\node[blacknode, scale=0.4] (c) at (-0.4,-2.6) {}; \node at (-0.9,-2.5) {$p=2$};
\node[whitenode, scale=0.4] (j) at (-0.6,-1.6) {};

\node at (-0.9,-1.3) {$B_2$};

\node[whitenode, scale=0.4] (i) at (0.31,-0.35) {};

\node at (0.6,-0.7) {$B_1$};

\node[blacknode, scale=0.4] (d) at (-0.17,-3.6) {}; \node at (0.05,-3.45) {$6$};
\node[whitenode, scale=0.4] (e) at (-1.25,-3.14) {};
\node[blacknode, scale=0.4] (f) at (-1.05,-4.1) {}; \node at (-1.3,-4) {$5$};

\node at (-1.5,-2.8) {$B_3$};

\node[whitenode, scale=0.4] (k) at (-1.05,-4.9) {};
\node[whitenode, scale=0.4] (l) at (-0.17,-4.34) {};

\node at (-1.6,-4.6) {$F=B_4$};
\node at (0.1,-4) {$B_5$};

\node[whitenode, scale=0.4] (g) at (1.2,-2) {};
\node[whitenode, scale=0.4] (h) at (1.5,-2.55) {}; \node at (1.3,-2.7) {$4$};

\node at (1.3,-1.7) {$B_6$};

\node[whitenode, scale=0.4] (m) at (1.5,-3.3) {};

\node at (1.8,-3) {$B_7$};

\node at (0,-5.25) {};

\draw[-, thick] (a) -- (b);
\draw[-, thick] (a) -- (c);
\draw[-, thick] (c) -- (b);
\draw[-, thick] (a) -- (j);
\draw[-, thick] (b) -- (j);
\draw[-, thick] (c) -- (j);

\draw[-, thick] (d) -- (e);
\draw[-, thick] (d) -- (f);
\draw[-, thick] (f) -- (e);
\draw[-, thick] (d) -- (c);
\draw[-, thick] (e) -- (c);
\draw[-, thick] (f) -- (c);

\draw[-, thick] (f) -- (k);
\draw[-, thick] (d) -- (l);

\draw[-, thick] (g) -- (h);
\draw[-, thick] (g) -- (b);
\draw[-, thick] (h) -- (b);

\draw[-, thick] (h) -- (m);

\draw[-, thick] (i) -- (a);

\draw[dashed, ultra thick] (-0.9,-2.2) .. controls (-0.2,-2.1) and (0.2,-2.5) .. (0.2,-3);

\end{tikzpicture}
\caption{}
\end{subfigure} \hspace{2mm}
\begin{subfigure}[h]{0.4\textwidth}
\centering

\begin{tikzpicture}[scale=1.5,
blacknode/.style={circle, draw=black!, fill=black, thick, minimum size= 8mm},
rednode/.style={circle, draw=black!, fill=red!80, thick, minimum size= 8mm},
whitenode/.style={circle, draw=black!, fill=white!, thick, minimum size= 8mm},
]

\node[whitenode, scale=0.4] (u_B1) at (1,-1) {}; \node at (1.4,-1) {$u_{B_1}$};
\node[whitenode, scale=0.4] (u_B2) at (0,-1) {}; \node at (-0.5,-1) {$u_{B_2}$};
\node[whitenode, scale=0.4] (u_B6) at (1,-2.6) {}; \node at (1.5,-2.6) {$u_{B_6}$};
\node[whitenode, scale=0.4] (u_B7) at (1.5,-4.4) {}; \node at (2,-4.4) {$u_{B_7}$};
\node[whitenode, scale=0.4] (u_B3) at (-0.5,-2.8) {}; \node at (-1,-2.8) {$u_{B_3}$};
\node[whitenode, scale=0.4] (u_B4) at (-1,-4.6) {}; \node at (-1.7,-4.6) {$u_F=u_{B_4}$};
\node[whitenode, scale=0.4] (u_B5) at (0,-4.6) {}; \node at (0.5,-4.6) {$u_{B_5}$};

\node[whitenode, scale=0.4] (v1) at (0.5,0) {}; \node at (-0.1,0) {$v_r=v_1$};
\node[whitenode, scale=0.4] (v2) at (-0.5,-1.8) {}; \node at (-1.1,-1.8) {$v_p=v_2$};
\node[whitenode, scale=0.4] (v3) at (0.5,-1.8) {}; \node at (1,-1.8) {$v_3$};
\node[whitenode, scale=0.4] (v4) at (1.5,-3.4) {}; \node at (2,-3.4) {$v_4$};
\node[whitenode, scale=0.4] (v5) at (-1,-3.6) {}; \node at (-1.5,-3.6) {$v_5$};
\node[whitenode, scale=0.4] (v6) at (0,-3.6) {}; \node at (0.5,-3.6) {$v_6$};

\draw[-, thick] (v1) -- (u_B1);
\draw[-, thick] (v1) -- (u_B2);
\draw[-, thick] (u_B2) -- (v2);
\draw[-, thick] (u_B2) -- (v3);
\draw[-, thick] (v3) -- (u_B6);
\draw[-, thick] (u_B6) -- (v4);
\draw[-, thick] (v4) -- (u_B7);
\draw[-, thick] (v2) -- (u_B3);
\draw[-, thick] (u_B3) -- (v5);
\draw[-, thick] (u_B3) -- (v6);
\draw[-, thick] (v5) -- (u_B4);
\draw[-, thick] (v6) -- (u_B5);

\end{tikzpicture}
\caption{}
\end{subfigure}
\caption{Figure (a) represents a twin-free block graph $G$ and Figure (b) represents $T_G$. The vertices underneath the dashed curve represent those deleted from $G$ to obtain $G'$. The black vertices represent vertices in the set $A$. (All notations are as in the proof of Theorem \ref{theorem_block}.)}
\label{block figure}
\end{figure}
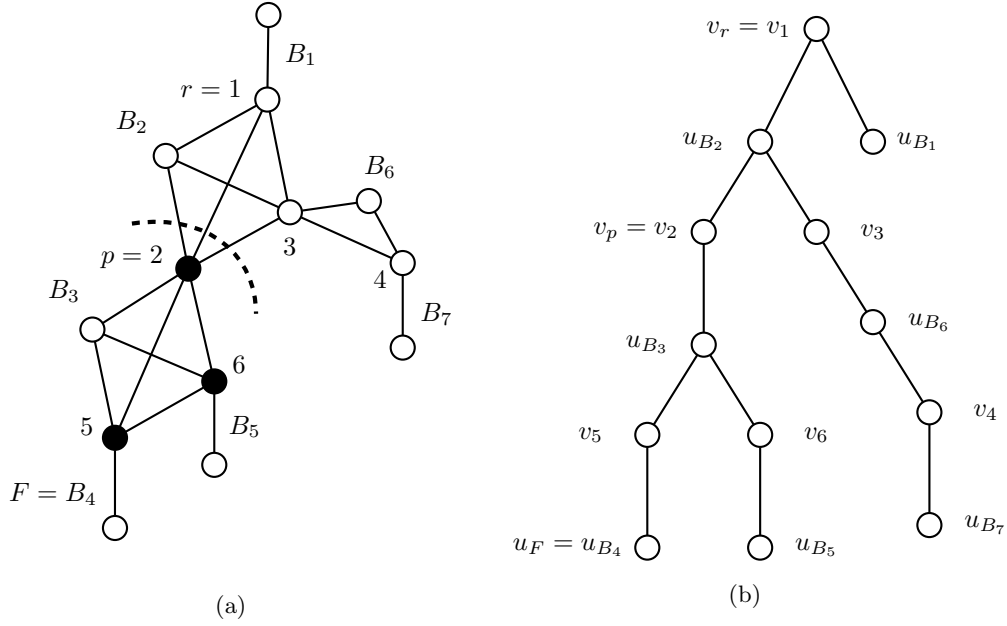


We now establish the following two claims related to the sets defined above.

\begin{claim} \label{claim_A LTD-set}
The set $A$ is an LTD-set of $G[U]$.
\end{claim}

\begin{claimproof}
That $A$ is a TD-set of $G[U]$ is clear from the structure of $G$. We show that $A$ is also a locating set of $G[U]$. So, let us assume that vertices $w,z \in U \setminus A$. This implies that both $w$ and $z$ are not cut-vertices of $G$. This further implies that $N_{G[U]}(w) \subset A$, or else, the block of $G$ that contains $w$ has another non-cut-vertex which is then a twin with $w$ in $G$, a contradiction. Similarly, $N_{G[U]}(z) \subset A$. In other words, $N_G(w) = N_{G[U]}(w) \cap A$ and $N_G(z) = N_{G[U]}(z) \cap A$. Thus, if $N_{G[U]}(w) \cap A = N_{G[U]} \cap A$, it implies that $N_G(w) = N_G(z)$ and thus, $w$ and $z$ are twins in $G$, again a contradiction. Hence, $A$ is a locating set of $G[U]$ and this proves the claim.
\end{claimproof}
\begin{claim} \label{claim_UA}
$|A| \leq \frac{2}{3}|U|$.
\end{claim}

\begin{claimproof}
Let $u_{B_1}, u_{B_2}, \ldots , u_{B_m}$ be $m \geq 1$ children of $v_p$ in $T_G$ and let each block $B_i$ of $G$ be of order $n_i$. Now, due to the twin-free nature of $G$, each vertex $u_{B_i}$ of $T_G$ has at least $n_i-2$ and at most $n_i-1$ children (and hence at least $n_i-2$ and at most $n_i-1$ grandchildren as well). To be more precise, assume that, for $0 \leq s \leq m$, the vertices $u_{B_1}, u_{B_2}, \ldots , u_{B_s}$ have exactly $n_1-2, n_2-2, \ldots , n_s-2$ children, respectively, in $T_G$; and that the vertices $u_{B_{s+1}}, u_{B_{s+2}}, \ldots , u_{B_m}$ have exactly $n_{s+1}-1, n_{s+2}-1, \ldots , n_m-1$ children, respectively, in $T_G$. This implies that we have the following equation.
$$|U|= 1-s + 2\sum_{1 \leq i \leq m} (n_i-1) = 1 -2m -s + 2\sum_{1 \leq i \leq m} n_i.$$

Moreover, we have
$$|A| = 1-s + \sum_{1 \leq i \leq m} (n_i-1) = 1 - m-s + \sum_{1 \leq i \leq m} n_i.$$

\noindent By combining the above two equations, therefore, we have
\begin{flalign*}
|U| - (2|A|-1) &= s \geq 0 \implies |A| \leq \frac{1}{2}(|U|+1) \leq \frac{2}{3}|U|,
\end{flalign*}
where the last inequality follows from noticing that $|U| \geq 3$. Hence, this proves the claim.
\end{claimproof}

Now, let $G' = G - U$, that is, $G'$ is the graph obtained by deleting from $G$ all vertices (and edges incident with them) in the blocks $B \in \mathcal{B}_p$. Notice that $G'$ is still a connected block graph; and assume that the order of $G'$ is $n'$ (which is strictly less than $n$). We next divide the proof according to whether $G'$ is twin-free, has twins or is isomorphic to a $3$-path.

\begin{case}[$G'$ is either twin-free or is isomorphic to a $3$-path]
We further subdivide this case into the following.

\begin{subcase}[$n' \leq 2$]
In this subcase, if $n'=2$, then the two vertices of $G'$ form an edge of $G'$ (since $G$ is connected). Hence, the two vertices of $G'$ are closed twins of degree $1$, contrary to our initial assumption in this case. Therefore, let us assume that $n' \leq 1$. If $n'=1$, then $v_p$ has no grandparent in $T_G$. In other words, there is no vertex of $v$-type above $v_p$ in $T_G$. This implies that $v_p$ must itself be the root vertex of $T_G$. However, this, in turn, implies that $G'$ is an empty graph which contradicts the fact that $n'=1$. Thus, we must have $n'=0$. In this case too, $v_p$ must itself be the root vertex of of $T_G$ and so, $G=G[U]$ and $|U|=n$. Therefore, by Claim~\ref{claim_A LTD-set}, the set $A$ is an LTD-set of $G[U]=G$. Moreover, by Claim~\ref{claim_UA}, we have
$$|A| \leq \frac{2}{3}|U| = \frac{2}{3}n.$$
\end{subcase}

\begin{subcase}[$n' \geq 3$] In this subcase, by the induction hypothesis, we have $\gamma^L_t(G') \leq \frac{2}{3}n'$. Suppose now that $S' \subset V(G')$ is a minimum LTD-set of $G'$, that is with $|S'|=\gamma^L_t(G')$. We then claim that the set $S = S' \cup A$ is an LTD-set of $G$. To prove so, we first see that the set $S$ is a TD-set of $G$, since $S'$ is a TD-set of $G'$ and $A$ is a TD-set of $G[U]$. Moreover, $S$ is also a locating set of $G$ due to the following two reasons.

\begin{enumerate}[(1)]
\item Any two distinct vertices $w \in V(G') \setminus S'$ and $z \in V(G) \setminus S$  are located by $S'$.
\item By Claim~\ref{claim_A LTD-set}, the set $A$ is a locating set of $G[U]$.
\end{enumerate}

Using Claim~\ref{claim_UA}, therefore, the two-thirds bound on $\gamma^L_t(G)$ in this subcase is established by the following inequality.
\begin{flalign*}
\gamma^L_t(G) \leq |S| &= |S'|+|A| \leq \gamma^L_t(G')+ \frac{2}{3} |U| \leq \frac{2}{3} \Big( n' + |U| \Big) = \frac{2}{3}n.
\end{flalign*}
\end{subcase}
\end{case}

We next turn to the case that $G'$ has twins.

\begin{case}[$G'$ is neither twin-free nor isomorphic to a $3$-path]
Assume that $x$ and $y$ are two vertices of $G'$ which are twins in $G'$. Then, without loss of generality, there exists an edge in $G$ between the vertices $p$ and $x$ and there is no edge in $G$ between $p$ and $y$. This implies that $x$ and $p$ belong to the same block $X$, say, of $G$ to which $y$ does not belong. Moreover, notice that there can be only one such $y$ that is a twin of $x$ in $G'$. Let $Y$ be a block of $G$ to which the vertex $y$ belongs. Next, we prove the following claim.

\begin{claim}\label{claim_C}
If $x$ does not belong to the block $Y$, then $\gamma^L_t(G) \le \frac{2}{3}n$.
\end{claim}

\begin{claimproof}
If $|V(Y)| \geq 3$, it would mean that $x$ would have at least two neighbors, say, $v$ and $w$ in $Y$ (since $x$ and $y$ are twins in $G'$). Then the vertices $v,w,x,y$ would induce a $K_4$ in the block graph $G$ making $x$ belong to $Y$, a contradiction. Therefore, we have $|V(Y)| = 2$, that is, $Y$ is a leaf block of $G$. So, let $V(Y) = \{v,y\}$. Therefore, we have $v \in N_{G'}(x)$. Now, assume that $deg_G(y) \geq 2$, that is, there exists some $w \in N_G(y) \setminus \{v\}$ (with $w = x$, possibly). Since $x$ and $y$ are twins in $G'$, if $w \neq x$, we have $v,w \in N_G(x)$ and thus, the set $\{v,w,x,y\}$ induces a $K_4$ in $G$. Moreover, if $w=x$, then the set $\{v,x,y\}$ induces a $K_3$ in $G$. Either way, the block $Y$ with $|V(Y)| = 2$ is contained inside the subgraph $K_3$ of $G$ which is a contradiction, since $Y$, being a block, is a maximal complete subgraph of $G$. Therefore, we have $deg_G(y) = 1$, that is, the vertex $y$ is a leaf with $v$ as its support vertex in $G$.

On the other hand, if $deg_G(x) \geq 3$, then $x$ must have a neighbor $w$, say, in $G$ other than $p$ and $v$. Therefore, $v,w \in N_{G'}(x)$. Since, $x$ and $y$ are twins in $G'$, we must also have $w \in N_{G}(y)$ making $deg_G(y) \geq 2$, a contradiction. Hence, we have $deg_G(x) = 2$ with $N_G(x) = \{p,v\}$. This implies that the set $\{x,v,y\}$ induces a $P_3$ in $G'$, where, the vertices $x$ and $y$ are leaves of $G'$ with $v$ as their common support vertex. Let $X_v$ be the block of $G$ containing the vertices $x$ and $v$. Notice that $X_v = X$, possibly, if $pv \in E(G)$ (recall that $X$ is the block of $G$ containing the vertices $x$ and $p$). Since $G' \not \cong P_3$ (by our assumption in this case), the vertex $v$ must belong to another block of $G$ other than $X_v$ and $Y$. Now, let $G'' = G' - x$. Thus,  $G''$ is still a connected block graph.  

 If the graph $G''$ had twins, then the vertex $v$ must be one of them. However, this possibility cannot arise as $y$ is a leaf of $G''$ with $v$ as its support vertex. Thus, $G''$ is also twin-free. Therefore, by our induction hypothesis, let $S'' \subset V(G'')$ be a minimum LTD-set of $G''$. Then, we have $|S''| = \gamma^L_t(G'') \leq \frac{2}{3}(n'-1)$. 
 
 We now claim that $S = A \cup S''$ is an LTD-set of $G$. Since, by Claim~\ref{claim_A LTD-set}, the set $A$ is an LTD-set of $G[U]$ and $S''$ is assumed to be an LTD-set of $G''$, to show that $S$ is an LTD-set of $G$, we have to show that
 \begin{enumerate}[(i)]
 \item the vertex $x$ is totally dominated and located from every other vertex of $V(G) \setminus S$ by $S$; and
 \item every pair $s,t$ with $s \in V(G'') \setminus S''$ and $t \in U \setminus A$ are located by $S$.
 \end{enumerate}
 To do so, we note that the support vertex $v$ must belong to the LTD-set $S''$ of $G''$. This implies that the vertex $x$ is totally dominated by $v \in S$ and is located from every vertex of $U \setminus A$ by $v \in S$. Moreover, the vertex $p \in A$. This implies that $x$ is located from every vertex of $V(G'') \setminus S''$ by $p \in S$. This implies that $x$ is located by $S$ from every vertex of $V(G) \setminus S$. Furthermore, let $s \in V(G'') \setminus S''$ and $t \in U \setminus A$. Since the vertex $t$ of $U \setminus A$ has a neighbor $t'$, say, that is a cut-vertex of $G$, it implies that $t' \in U$. Hence, we have $t' \in A$. This implies that $t' \in S$ locates the pair $s,t$. This proves that $S$ is an LTD-set of $G$. Hence, by Claim~\ref{claim_UA}, we have
\begin{flalign*}
\gamma^L_t(G) \leq |S| = |S''|+|A| \leq \frac{2}{3} \Big( (n'-1) + |U| \Big) < \frac{2}{3}n.
\end{flalign*}
This proves the claim.
\end{claimproof}

In view of Claim~\ref{claim_C} therefore, we assume for the rest of this proof that $x$ also belongs to the block $Y$ of $G$. Now, clearly, $X \neq Y$, or else, $py$ would be an edge in $G$, contradicting our earlier observation. Thus, $x$ is a cut-vertex of $G$ belonging to the distinct blocks $X$ and $Y$ of $G$. Moreover, $|V(X)|=2$, or else, again, $x$ and $y$ would not be twins in $G'$, a contradiction. More precisely, $V(X)=\{x,p\}$. We also observe here that $y$ cannot be a cut-vertex of $G$, or else, $x$ and $y$ would not be twins in $G'$, again the contradiction as before. Therefore, since $G$ is twin-free, every vertex other than $y$ of the block $Y$ must be a cut-vertex of~$G$.

Now, we again look at the block graph $G'' = G'- x$ on, say, $n'' ~(=n'-1)$ vertices
(the graph induced by the vertices on the right of the dashed curve in Figure~\ref{Figure 2a}). Notice that, in the tree $T_G$, the vertex $v_x$ cannot have any children other than $u_X$, or else, $x$ and $y$ cannot be twins in $G'$, contrary to our assumption for this case. This implies that the block graph $G''$ is also connected. We also have $y \in V(G'')$. Thus, $n'' \geq 1$. However, we can show that $n'' \neq 2$. Suppose, to the contrary, that  $n'' = 2$. In this case, the graph $G'$ is $K_3$ which implies that $G$ has twins, a contradiction. Next, we divide this case into the following subcases according to the order $n''$ of $G''$.


\begin{figure}[t!]
\centering

\begin{subfigure}[h]{0.4\textwidth}
\centering

\begin{tikzpicture}[scale=1.5,
blacknode/.style={circle, draw=black!, fill=black, thick, minimum size= 8mm},
rednode/.style={circle, draw=black!, fill=red!80, thick, minimum size= 8mm},
whitenode/.style={circle, draw=black!, fill=white!, thick, minimum size= 8mm},
]

\node[whitenode, scale=0.4] (b) at (0.5,-2.1) {}; \node at (0.5,-2.4) {$x$};
\node[blacknode, scale=0.4] (c) at (-0.4,-2.6) {}; \node at (-0.7,-2.5) {$p$};


\node[blacknode, scale=0.4] (d) at (-0.17,-3.6) {}; 
\node[whitenode, scale=0.4] (e) at (-1.25,-3.14) {};
\node[blacknode, scale=0.4] (f) at (-1.05,-4.1) {};

\node[whitenode, scale=0.4] (k) at (-1.05,-4.9) {};
\node[whitenode, scale=0.4] (l) at (-0.17,-4.34) {};

\node[whitenode, fill=black!50, scale=0.4] (g) at (1.2,-2) {}; \node at (1.4,-1.8) {$y$};
\node[whitenode, fill=black!50, scale=0.4] (h) at (1.5,-2.55) {}; \node at (1.3,-2.7) {$z$};

\node[whitenode, scale=0.4] (m) at (1.5,-3.3) {}; \node at (1.8,-3.3) {$x''$};

\node at (0,-5.25) {};

\draw[-, thick] (c) -- (b);

\draw[-, thick] (d) -- (e);
\draw[-, thick] (d) -- (f);
\draw[-, thick] (f) -- (e);
\draw[-, thick] (d) -- (c);
\draw[-, thick] (e) -- (c);
\draw[-, thick] (f) -- (c);

\draw[-, thick] (f) -- (k);
\draw[-, thick] (d) -- (l);

\draw[-, thick] (g) -- (h);
\draw[-, thick] (g) -- (b);
\draw[-, thick] (h) -- (b);

\draw[-, thick] (h) -- (m);

\draw[dashed, ultra thick] (0,-2) .. controls (1,-1.3) and (1,-2) .. (0.5,-2.8);

\draw[-, thick, dotted] (-0.9,-2.3) -- (0.5,-1.6);
\draw[-, thick, dotted] (-0.4,-3.2) -- (1,-2.5);
\draw[-, thick, dotted] (-0.9,-2.3) -- (-0.4,-3.2);
\draw[-, thick, dotted] (0.5,-1.6) -- (1,-2.5);

\node at (-0.5,-1.8) {$X$};

\draw[-, thick, dotted] (0.2,-1.5) -- (1.8,-1.5);
\draw[-, thick, dotted] (0.2,-2.9) -- (1.8,-2.9);
\draw[-, thick, dotted] (0.2,-1.5) -- (0.2,-2.9);
\draw[-, thick, dotted] (1.8,-1.5) -- (1.8,-2.9);

\node at (2.1,-2.2) {$Y$};

\draw[-, thick, dotted] (1.2,-2.3) -- (1.9,-2.3);
\draw[-, thick, dotted] (1.2,-3.6) -- (1.9,-3.6);
\draw[-, thick, dotted] (1.2,-2.3) -- (1.2,-3.6);
\draw[-, thick, dotted] (1.9,-3.6) -- (1.9,-2.3);

\node at (1.6,-3.8) {$X''$};

\end{tikzpicture}
\caption{The vertices to the left of the dashed curve represent those deleted from $G$ to obtain $G''$. $G'' \cong P_3$. The black vertices constitute the set $A$ and the grey vertices constitute an LTD-set $S''$ of $G''$.\\} \label{Figure 2a}
\end{subfigure} \hspace{2mm}
\begin{subfigure}[h]{0.4\textwidth}
\centering

\begin{tikzpicture}[scale=1.5,
blacknode/.style={circle, draw=black!, fill=black, thick, minimum size= 8mm},
rednode/.style={circle, draw=black!, fill=red!80, thick, minimum size= 8mm},
whitenode/.style={circle, draw=black!, fill=white!, thick, minimum size= 8mm},
]

\node[blacknode, scale=0.4] (b) at (0.5,-2.1) {}; \node at (0.5,-2.4) {$x$};
\node[blacknode, scale=0.4] (c) at (-0.4,-2.6) {}; \node at (-0.7,-2.5) {$p$};


\node[blacknode, scale=0.4] (d) at (-0.17,-3.6) {}; 
\node[whitenode, scale=0.4] (e) at (-1.25,-3.14) {};
\node[blacknode, scale=0.4] (f) at (-1.05,-4.1) {};

\node[whitenode, scale=0.4] (k) at (-1.05,-4.9) {};
\node[whitenode, scale=0.4] (l) at (-0.17,-4.34) {};

\node[whitenode, scale=0.4] (g) at (1.2,-2) {}; \node at (1.4,-1.8) {$y$};
\node[whitenode,  fill=black!50, scale=0.4] (h) at (1.5,-2.55) {}; \node at (1.3,-2.7) {$z$};

\node[whitenode, scale=0.4] (m) at (1.5,-3.3) {}; \node at (1.8,-3.3) {$x''$};

\node[whitenode,  fill=black!50, scale=0.4] (n) at (2.5,-2.55) {};
\node[whitenode, scale=0.4] (o) at (2,-2.9) {};

\node[whitenode, scale=0.4] (p) at (2.5,-3.3) {};

\node at (0,-5.25) {};

\draw[-, thick] (c) -- (b);

\draw[-, thick] (d) -- (e);
\draw[-, thick] (d) -- (f);
\draw[-, thick] (f) -- (e);
\draw[-, thick] (d) -- (c);
\draw[-, thick] (e) -- (c);
\draw[-, thick] (f) -- (c);

\draw[-, thick] (f) -- (k);
\draw[-, thick] (d) -- (l);

\draw[-, thick] (g) -- (h);
\draw[-, thick] (g) -- (b);
\draw[-, thick] (h) -- (b);

\draw[-, thick] (h) -- (m);

\draw[-, thick] (h) -- (n);
\draw[-, thick] (h) -- (o);
\draw[-, thick] (n) -- (o);

\draw[-, thick] (n) -- (p);



\draw[dashed, ultra thick] (2.5,-1.8) .. controls (1,-2) and (1,-2.8) .. (1,-3.5);

\draw[-, thick, dotted] (-0.9,-2.3) -- (0.5,-1.6);
\draw[-, thick, dotted] (-0.4,-3.2) -- (1,-2.5);
\draw[-, thick, dotted] (-0.9,-2.3) -- (-0.4,-3.2);
\draw[-, thick, dotted] (0.5,-1.6) -- (1,-2.5);

\node at (-0.5,-1.8) {$X$};

\draw[-, thick, dotted] (0.2,-1.5) -- (1.8,-1.5);
\draw[-, thick, dotted] (0.2,-2.9) -- (1.8,-2.9);
\draw[-, thick, dotted] (0.2,-1.5) -- (0.2,-2.9);
\draw[-, thick, dotted] (1.8,-1.5) -- (1.8,-2.9);

\node at (2.1,-2.2) {$Y$};

\draw[-, thick, dotted] (1.2,-2.3) -- (1.9,-2.3);
\draw[-, thick, dotted] (1.2,-3.6) -- (1.9,-3.6);
\draw[-, thick, dotted] (1.2,-2.3) -- (1.2,-3.6);
\draw[-, thick, dotted] (1.9,-3.6) -- (1.9,-2.3);

\node at (1.6,-3.8) {$X''$};

\end{tikzpicture}
\caption{The vertices to the left of the dashed curve represent those deleted from $G$ to obtain $G^\star$. $G''$ has twins $x''$ and $y$; and $G^\star$ is a twin-free block graph. The black vertices constitute the set $A \cup \{x\}$ and the grey vertices constitute an LTD-set $S^\star$ of $G^\star$.} \label{Figure 2b}
\end{subfigure}

\caption{Twin-free block graph $G$. The dotted boxes mark the blocks $X$, $X''$ and $Y$ of $G$ as in the proof of Theorem \ref{theorem_block}.}
\label{block figure 2}
\end{figure}
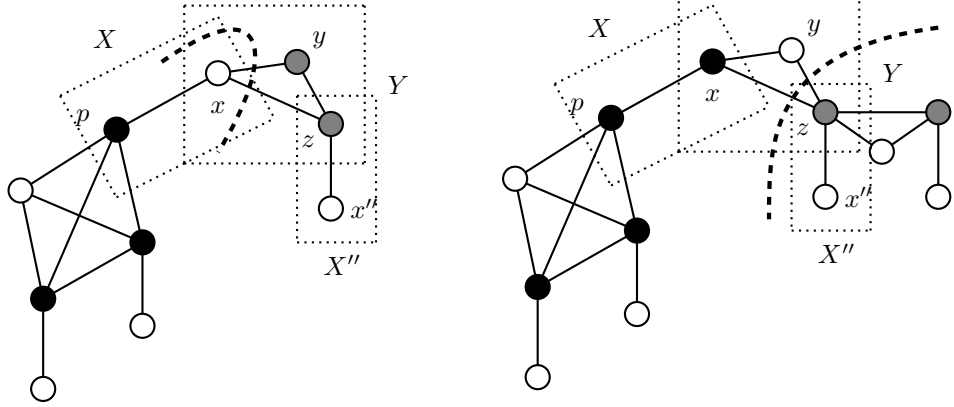


\begin{subcase}[$n'' = 1$]
In this subcase, we claim that $S=A \cup \{x\}$ is an LTD-set of $G$. It is clear that $S$ is a TD-set of $G$; by Claim~\ref{claim_A LTD-set}, set $A$ is an LTD-set of $G[U]$. The vertex $y$ is located from any vertex in $U \setminus A$ by the vertex $x$, and thus the set $S$ is also locating. Therefore, in this case, using Claim~\ref{claim_UA} we have
\begin{flalign*}
\gamma^L_t(G) \leq |S| = |A|+1 < \frac{2}{3} \Big( |U| + 2 \Big) = \frac{2}{3}n.
\end{flalign*}
\end{subcase}

\begin{subcase}[$n'' \geq 3$ and $G''$ is either twin-free or is isomorphic to a $3$-path]
Since $n''$ is at least $3$ and is strictly less than $n$, by the induction hypothesis, we have $\gamma^L_t(G'') \leq \frac{2}{3} n''$. Moreover, let $S''$ be a minimum LTD-set of $G''$, that is with $|S''|=\gamma^L_t(G'')$. We next claim the following.

\begin{claim}
The set $S=S'' \cup A$ is an LTD-set of $G$.
\end{claim}

\begin{claimproof}
Since $S''$ is a TD-set of $G''$ and $A$ is a TD-set of $G[U \cup \{x\}]$, the set $S$ is therefore a TD-set of $G$. Next we show that $S$ is also a locating set of $G$. To begin with, we note that $Y''=Y-x$ is a block of $G''$ containing the vertex $y$. Now, since $y$ is not a cut-vertex of $G$, we have $S'' \cap Y'' \neq \emptyset$ (or else, $y$ is not dominated by $S''$). This implies that $x$ is located by $S''$ from all vertices in $U \setminus A$. Moreover, $x$ is also located by $p$ from all vertices in $V(G'') \setminus S''$. Next, any pair $w,z$ of distinct vertices with $w \in V(G'') \setminus S''$ and $z \in V(G) \setminus (S \cup \{x\})$ are located by $S''$. Finally, any distinct pair of vertices $w,z \in U \setminus A$ are located by $A$, since the latter is an LTD-set of $G[U]$ by Claim~\ref{claim_A LTD-set}.
\end{claimproof}

Therefore, in this subcase, once again by Claim~\ref{claim_UA}, the theorem follows from the following inequality:
\begin{flalign*}
\gamma^L_t(G) \leq |S| &= |S''|+|A| = \gamma^L_t(G'')+ |A| \leq \frac{2}{3} \Big( n'' + |U| \Big) < \frac{2}{3}n.
\end{flalign*}
\end{subcase}

\begin{subcase}[$n'' \geq 3$ and $G''$ is neither twin-free nor is isomorphic to a $3$-path]
Assume that $x''$ and $y''$ are a pair of twins of $G''$. Moreover, for $x''$ and $y''$ to be twins in $G''$, at least one of them must be in the block $Y$. Let us, without loss of generality, assume that $y'' \in V(Y)$. 

We next observe that the vertices $y$ and $y''$ are the same. To prove so, suppose, to the contrary, that $y'' \neq y$. Then $y''$ is a cut-vertex of $G$ and so, for $x''$ and $y''$ to be twins in $G''$, $x''$ must \emph{not} belong to the block $Y$ of $G$. However, this, in turn, implies that $y$ is a neighbor of $y''$ but not of $x''$ and so, $x''$ and $y''$ are not twins in $G''$, a contradiction all the same. This, therefore, proves the observation.

Again, the vertex $x'' \notin Y$, since otherwise, $x'' \neq y'' = y$ implies that $x''$ is a cut-vertex of $G$, thus forcing $x''$ and $y''$ to not be twins, contrary to our supposition. Let $x''$ belong to the block $X'' ~(\neq Y)$ of $G''$ (and of $G$). We now try to establish the structure of the block $Y$ of $G$. Notice that, by the structure of a block graph, the twins $x''$ and $y$ in $G''$ must have a \emph{single} common neighbor $z$, say, in $G''$ such that $z$ is a cut-vertex of $G$ belonging to both the blocks $Y$ and $X''$ of $G$. Furthermore, if the block $Y$ contains any vertex of $G$ other than the vertices $x,y$ and $z$, then $x''$ and $y$ are not twins in $G''$, a contradiction. Thus, we have $V(Y) = \{x,y,z\}$.

Next, to understand the structure of the block $X''$ of $G''$, we see that neither can $X''$ contain any vertex other than $z$ and $x''$, nor can $x''$ be a cut-vertex of $G$; or else, we again have the contradiction that $x''$ and $y$ are not twins in $G''$. Therefore, this implies that $V(X'') = \{x'',z\}$, that is, $X''$ is a leaf block of $G''$ (and of $G$). See Figure~\ref{block figure 2} for the structure of the blocks $X''$ and $Y$.

With that, we look at the block graph $G^\star = G''-y$ (the graph induced by the vertices on the right of the dashed curve in Figure~\ref{Figure 2b}). Then, $G^\star$ is again a connected graph, since $y$ is not a cut-vertex of $G$. Moreover, the order $n^\star$ of $G^\star$ is at least~2 (since $x'',z \in V(G^\star)$). If, however, $n^\star=2$, then we have $V(G'') = \{x'',y,z\}$ and thus, $G''$ is isomorphic to a $3$-path, contrary to our assumption in this subcase. Therefore, we have $n^\star \geq 3$. We next show the following claim.

\begin{claim}
The graph $G^\star$ is twin-free.
\end{claim}

\begin{claimproof}
Suppose, to the contrary, that the block graph $G^\star$ has a pair of twins. In this case, one of them must be the cut-vertex $z$ of $G$. Let $x^\star$ be the other vertex of $G^\star$ such that $x^\star$ and $z$ are twins in $G^\star$. Since $x''$ is a neighbor of $z$ alone in $G^\star$, therefore $z$ cannot be a twin in $G^\star$ of any vertex other than $x''$. In other words, $x^\star=x''$. However, since $deg_{G^\star} (x'') = 1$, we have $deg_{G^\star} (z) = 1$ and, hence, the graph $G^\star$ is simply the edge $x''z$ of $G$. This however, contradicts the fact that $n^\star \geq 3$. Hence, this proves that $G^\star$ is twin-free.
\end{claimproof}

\noindent Since $n^\star$ is at least~3 and is strictly less than $n$, by the induction hypothesis, we have $\gamma^L_t(G^\star) \leq \frac{2}{3} n^\star$. Moreover, let $S^\star$ be a minimum LTD-set of $G^\star$, that is with $|S^\star|=\gamma^L_t(G^\star)$. We next claim the following.

\begin{claim}
The set $S=S^\star \cup A \cup \{x\}$ is an LTD-set of $G$.
\end{claim}

\begin{claimproof}
Since $S^\star$ is a TD-set of $G^\star$ and $A \cup \{x\}$ is a TD-set of $G[U \cup \{x,y\}]$, the set $S$ is therefore a TD-set of $G$. Next we show that $S$ is also a locating set of $G$. To begin with, we note that, since $x''$ is a leaf in $G^\star$, its support vertex $z$ must be in the LTD-set $S^\star$ of $G^\star$. Thus, the vertex $y$ is located from every other vertex in $V(G) \setminus S$ by the set $\{x,z\}$. Next, any pair $w_1,w_2$
of distinct vertices with $w_1 \in V(G^\star) \setminus S^\star$ and $w_2 \in V(G) \setminus S$, respectively, are located by $S^\star$. Finally, by Claim~\ref{claim_A LTD-set}, any pair of distinct vertices $w_1,w_2 \in U \setminus A$ are located by the set $A$.
\end{claimproof}

\noindent Therefore, again using Claim~\ref{claim_UA}, in this subcase, the theorem follows from the following inequality:
\begin{flalign*}
\gamma^L_t(G) \leq |S| &= |S^\star|+|A|+1 < \frac{2}{3} \Big( n^\star + |U|+ 2 \Big) = \frac{2}{3}n.
\end{flalign*}

\end{subcase}
\end{case}
This completes the proof.
\end{proof}

\begin{figure}[htb]
\begin{center}
\begin{tikzpicture}[scale=.7,style=thick,x=1cm,y=1cm]
\def\vr{5pt}
\path (0,0) coordinate (v1);
\path (0,1) coordinate (v2);
\path (0,2) coordinate (v3);
\draw (v1) -- (v2);
\draw (v2) -- (v3);
\path (2,0) coordinate (u1);
\path (2,1) coordinate (u2);
\path (2,2) coordinate (u3);
\draw (u1) -- (u2);
\draw (u2) -- (u3);
\path (4,0) coordinate (w1);
\path (4,1) coordinate (w2);
\path (4,2) coordinate (w3);
\draw (w1) -- (w2);
\draw (w2) -- (w3);
\path (6,0) coordinate (x1);
\path (6,1) coordinate (x2);
\path (6,2) coordinate (x3);
\draw (x1) -- (x2);
\draw (x2) -- (x3);
\path (8,0) coordinate (y1);
\path (8,1) coordinate (y2);
\path (8,2) coordinate (y3);
\draw (y1) -- (y2);
\draw (y2) -- (y3);
\path (10,0) coordinate (z1);
\path (10,1) coordinate (z2);
\path (10,2) coordinate (z3);
\draw (z1) -- (z2);
\draw (z2) -- (z3);
\draw (v3) -- (u3);
\draw (u3) -- (w3);
\draw (w3) -- (x3);
\draw (x3) -- (y3);
\draw (y3) -- (z3);
\draw (v1) [fill=white] circle (\vr);
\draw (v2) [fill=black] circle (\vr);
\draw (v3) [fill=black] circle (\vr);
\draw (u1) [fill=white] circle (\vr);
\draw (u2) [fill=black] circle (\vr);
\draw (u3) [fill=black] circle (\vr);
\draw (w1) [fill=white] circle (\vr);
\draw (w2) [fill=black] circle (\vr);
\draw (w3) [fill=black] circle (\vr);
\draw (x1) [fill=white] circle (\vr);
\draw (x2) [fill=black] circle (\vr);
\draw (x3) [fill=black] circle (\vr);
\draw (y1) [fill=white] circle (\vr);
\draw (y2) [fill=black] circle (\vr);
\draw (y3) [fill=black] circle (\vr);
\draw (z1) [fill=white] circle (\vr);
\draw (z2) [fill=black] circle (\vr);
\draw (z3) [fill=black] circle (\vr);
\draw (v3) to[out=30,in=150, distance=0.75cm ] (w3);
\draw (v3) to[out=40,in=100, distance=1cm ] (x3);
\draw (v3) to[out=50,in=110, distance=1.25cm ] (y3);
\draw (v3) to[out=60,in=120, distance=1.5cm ] (z3);
\draw (u3) to[out=30,in=150, distance=0.75cm ] (x3);
\draw (u3) to[out=40,in=100, distance=1cm ] (y3);
\draw (u3) to[out=50,in=110, distance=1.25cm ] (z3);
\draw (w3) to[out=30,in=150, distance=0.75cm ] (y3);
\draw (w3) to[out=40,in=100, distance=1cm ] (z3);
\draw (x3) to[out=30,in=150, distance=0.75cm ] (z3);

\end{tikzpicture}
\end{center}
\begin{center}
\caption{The $2$-corona $K_6 \circ P_2$ of a complete graph of order $6$.}
\label{2corona-K6}
\end{center}
\end{figure}
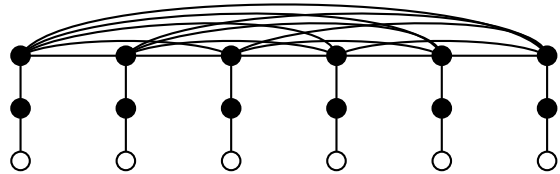

The ``twin-free'' condition for block graphs is necessary as, without it, the conjecture does not hold: for example, for $\Delta$-stars $K_{1,\Delta}$ with $\Delta \geq 3$ (which are block graphs), the locating-total domination number is $\Delta$. On the other hand, for any block graph $H$ of order $k\geq 2$, the $2$-corona $G = H \circ P_2$ is a twin-free block graph of order~$n = 3k$ and by Observation~\ref{obs:conj-tight}, it has locating-total domination number equal to its total domination number, that is, $\gamma_t^L(G) = \gamma_t(G) = 2k = \frac{2}{3}n$. See Figure~\ref{2corona-K6} for an illustration with $H$ a complete graph. Thus, we obtain the following.

\begin{proposition}
\label{tightness-of-block}
There are infinitely many connected twin-free block graphs $G$ of order $n$ with $\gamma_t^L(G) = \frac{2}{3}n$.
\end{proposition}

\section{Subcubic graphs}\label{sec:subcubic}

In this section, we establish a tight upper bound on the locating-total domination number of a subcubic graph, where a subcubic graph is a graph with maximum degree at most~$3$. For this purpose, let $\cF_\tdom$ be the family consisting of the three complete graphs $K_1$, $K_2$, and $K_4$, and a star $K_{1,3}$, that is,
\[
\cF_\tdom = \{K_1,K_2,K_4,K_{1,3}\}.
\]

Recall that a \emph{diamond} is the graph $K_4 - e$ where $e$ is an arbitrary edge of the $K_4$. A \emph{paw} is the graph obtained from a triangle $K_3$ by adding a new vertex and joining it with an edge to one vertex of the triangle. Equivalently, a paw is obtained from $K_{1,3}$ by adding an edge between two leaves.

For $k \ge 2$, we say a graph $G$ contains a $(d_1,d_2,\ldots,d_k)$-\emph{sequence} if there exists a path $v_1v_2 \ldots v_k$ such that $\deg_G(v_i) = d_i$ for all $i \in [k]$. We are now in a position to prove the following upper bound on the locating-total domination number of a subcubic graph.

\begin{theorem}
\label{t:thm-subcubic}
If $G \notin \cF_\tdom$ is a connected subcubic graph of order~$n \ge 3$, then $\gamma_t^L(G) \le \frac{2}{3}n$.
\end{theorem}
\begin{proof} Suppose, to the contrary, that the theorem is false. Among all counterexamples, let $G$ be one of minimum order~$n$. If $n = 3$, then $G \cong P_3$ or $G \cong K_3$, and in both cases $\gamma_t^L(G) = 2 = \frac{2}{3}n$, a contradiction. Hence, $n \ge 4$. Suppose $n = 4$. By assumption, $G \notin \{K_4,K_{1,3}\}$. If $G$ is a diamond or a paw, then let $S$ consist of one vertex of degree~$2$ and one vertex of degree~$3$, and if $G$ is a path or a cycle, then let $S$ consist of two adjacent vertices of degree~$2$. In all cases, $S$ is an LTD-set of $G$ of cardinality~$2$, and so $\gamma_t^L(G) \le 2 < \frac{2}{3}n$, a contradiction.  Hence, $n \ge 5$.

Suppose that $n = 5$. If $G$ is a path $P_5$ or a cycle $C_5$, then $\gamma_t^L(G) = 3 < \frac{2}{3}n$ (choose three consecutive vertices of degree~2), a contradiction. Hence, $\Delta(G) = 3$. Let $v$ be a vertex of degree~$3$ in $G$ with neighbors $v_1,v_2,v_3$. Let $v_4$ be the remaining vertex of $G$. Since $G$ is connected, we may assume, renaming vertices if necessary, that $v_1v_4$ is an edge. The set $\{v,v_1,v_2\}$ is an LTD-set of $G$, and so $\gamma_t^L(G) \le 3 < \frac{2}{3}n$, a contradiction.  Hence, $n \ge 6$.

Suppose that $n = 6$. If $G$ is a path $P_6$ or a cycle $C_6$, then $\gamma_t^L(G) = 4 = \frac{2}{3}n$ (choose four consecutive vertices of degree~2), a contradiction. Hence, $\Delta(G) = 3$. Let $v$ be a vertex of degree~$3$ in $G$ with neighbors $v_1,v_2,v_3$, and let $v_4$ and $v_5$ be the two remaining vertices of $G$. Since $G$ is connected, we may assume, renaming vertices if necessary, that $v_1v_4$ is an edge.
One of the sets $\{v,v_1,v_2,v_4\}$ and $\{v,v_1,v_3,v_4\}$ is an LTD-set of $G$, and so $\gamma_t^L(G) \le 4 = \frac{2}{3}n$, a contradiction. Hence, $n \ge 7$.

In what follows, we adopt the notation that if there is a $(d_1,d_2,\ldots,d_k)$-sequence in $G$, then $P \colon v_1v_2 \ldots v_k$ denotes a path in $G$ associated with such a sequence, where $\deg_G(v_i) = d_i$ for all $i \in [k]$. Further, we let $G' = G - V(P)$ and let $G'$ have order~$n'$, and so $n' = n - k$. Recall that $n \ge 7$.

We show firstly that there is no vertex of degree~$1$.

\begin{claim}
\label{c:claim1}
$\delta(G) \ge 2$.
\end{claim}
\begin{claimproof} Suppose, to the contrary, that $\delta(G) = 1$. We proceed further with a series of structural properties of the graph $G$ that show that certain $(d_1,d_2,\ldots,d_k)$-sequences are forbidden.

\begin{subclaim}
\label{c:claim1.1}
The following properties hold in the graph $G$.
\begin{enumerate}
\item[{\rm (a)}] There is no $(1,3,1)$-sequence.
\item[{\rm (b)}] There is no $(1,2,2)$-sequence.
\item[{\rm (c)}] There is no $(1,2,3,1)$-sequence.
\item[{\rm (d)}] There is no $(1,2,3,2,1)$-sequence.
\item[{\rm (e)}] There is no $(1,2,3)$-sequence.
\item[{\rm (f)}] There is no $(1,2)$-sequence.
\item[{\rm (g)}] There is no $(1,3,2)$-sequence.
\item[{\rm (h)}] There is no $(1,3,3,1)$-sequence.
\end{enumerate}
\end{subclaim}
\begin{subclaimproof} (a) Suppose that there is a $(1,3,1)$-sequence in $G$. In this case, $n' = n - 3 \ge 4$. Since $G$ is connected, so too is the graph $G'$. Let $v'$ be the third neighbor of $v_2$ in $G$ not on the path $P$. Suppose $G' \in \cF_\tdom$, implying that $G' \cong K_{1,3}$ with $v'$ as a leaf in $G'$. The graph $G$ is therefore determined, and has order~$n=7$. In this case, choosing $S$ to consist of the two support vertices (of degree~$3$) and a leaf neighbor of each support vertex produces an LTD-set of $G$ of cardinality~$4$, and so $\gamma_t^L(G) \le 4 < \frac{2}{3}n$. Hence, $G' \notin \cF_\tdom$. Since $G'$ is not a counterexample, it holds that $\gamma_t^L(G') \le \frac{2}{3}n' = \frac{2}{3}n - 2$. Every $\gamma_t^L$-set of $G'$ can be extended to an LTD-set of $G$ by adding to it the vertices $v_2$ and $v_3$, implying that $\gamma_t^L(G) \le \gamma_t^L(G') + 2 \le \frac{2}{3}n$, a contradiction.

(b) Suppose that there is a $(1,2,2)$-sequence in $G$. Let $v'$ be the second neighbor of $v_3$. As in the previous case, $n' = n - 3 \ge 4$ and $G'$ is connected. By part~(a), there is no $(1,3,1)$-sequence, implying that $G' \notin \cF_\tdom$ and $\gamma_t^L(G') \le \frac{2}{3}n' = \frac{2}{3}n - 2$. As before every $\gamma_t^L$-set of $G'$ can be extended to an LTD-set of $G$ by adding to it the vertices $v_2$ and $v_3$, implying that $\gamma_t^L(G) \le \frac{2}{3}n$, a contradiction.

(c) Suppose that there is a $(1,2,3,1)$-sequence in $G$. In this case, $n' = n - 4 \ge 3$ and $G'$ is connected. By part~(a), there is no $(1,3,1)$-sequence, implying that $G' \notin \cF_\tdom$ and $\gamma_t^L(G') \le \frac{2}{3}n' = \frac{2}{3}(n-4) < \frac{2}{3}n - 2$. Every $\gamma_t^L$-set of $G'$ can be extended to an LTD-set of $G$ by adding to it the vertices $v_2$ and $v_3$, implying that $\gamma_t^L(G) \le \frac{2}{3}n$, a contradiction.

(d) Suppose that there is a $(1,2,3,2,1)$-sequence in $G$. In this case, $G'$ is connected and $n' = n - 5 \ge 2$.
If $G' \in \cF_\tdom$, then $G' \cong K_2$ by the fact that there is no $(1,3,1)$-sequence in $G$ by part (a). The graph $G$ is therefore determined, and is obtained from a star $K_{1,3}$ by subdividing every edge once. We note that $G$ has order~$n=7$ and the set $N[v_3]$ (of non-leaves of $G$) is an LTD-set of $G$, implying that $\gamma_t^L(G) \le 4 < \frac{2}{3}n$, a contradiction. Hence, $G' \notin \cF_\tdom$. Thus, $\gamma_t^L(G') \le \frac{2}{3}n' = \frac{2}{3}(n-5) < \frac{2}{3}n - 3$. Every $\gamma_t^L$-set of $G'$ can be extended to an LTD-set of $G$ by adding to it the vertices $v_2, v_3$, and $v_4$, implying that $\gamma_t^L(G) < \frac{2}{3}n$, a contradiction.

(e) Suppose that there is a $(1,2,3)$-sequence in $G$. In this case, $n' = n - 3 \ge 4$ and $G'$ contains at most two components. Let $v_4$ and $v_4'$ be the two neighbors of $v_3$ different from $v_2$. By our earlier observations, each of $v_4$ and $v_4'$ has degree at least~$2$ in $G$, and therefore degree at least~$1$ in $G'$.

Suppose that $G'$ is disconnected. In this case, since there is no $(1,3,1)$-sequence, no $(1,2,3,1)$-sequence, and no $(1,2,3,2,1)$-sequence in $G$, neither component of $G'$ belongs to $\cF_\tdom$. By linearity, we therefore have that $\gamma_t^L(G') \le \frac{2}{3}n' = \frac{2}{3}n - 2$. Every $\gamma_t^L$-set of $G'$ can be extended to an LTD-set of $G$ by adding to it the vertices $v_2$ and $v_3$, implying that $\gamma_t^L(G) \le \frac{2}{3}n$, a contradiction. Hence, $G'$ is connected. Recall that $n' \ge 4$.

Suppose now that $G'$ is connected. If $G' \in \cF_\tdom$, then $G' \cong K_{1,3}$. Let $v_5$ be the central vertex of $G'$, and so each of $v_4$ and $v_4'$ is a leaf neighbor of $v_5$ in $G'$. The graph $G$ is therefore determined and $n = 7$. The set $\{v_2,v_3,v_4,v_5\}$ is an LTD-set of $G$, implying that $\gamma_t^L(G) \le 4 < \frac{2}{3}n$, a contradiction. Hence, $G' \notin \cF_\tdom$. Thus, $\gamma_t^L(G') \le \frac{2}{3}n' = \frac{2}{3}n - 2$. Every $\gamma_t^L$-set of $G'$ can be extended to an LTD-set of $G$ by adding to it the vertices $v_2$ and $v_3$, implying that $\gamma_t^L(G) \le \frac{2}{3}n$, a contradiction.

(f) Since there is no $(1,2,1)$-sequence (since $n\ge 7$), no $(1,2,2)$-sequence by (b) and no $(1,2,3)$-sequence by (e), there can be no $(1,2)$-sequence in $G$. Hence, part~(f) follows immediately from parts~(b) and~(e).

(g) Suppose that there is a $(1,3,2)$-sequence in $G$. In this case, $n' = n - 3 \ge 4$. If $G'$ is disconnected, then by parts~(a)--(f), neither component of $G'$ belongs to $\cF_\tdom$. By linearity, we therefore have that $\gamma_t^L(G') \le \frac{2}{3}n' = \frac{2}{3}n - 2$. Every $\gamma_t^L$-set of $G'$ can be extended to an LTD-set of $G$ by adding to it the vertices $v_2$ and $v_3$, implying that $\gamma_t^L(G) \le \frac{2}{3}n$, a contradiction. Hence, $G'$ is connected. If $G' \in \cF_\tdom$, then $G' \cong K_{1,3}$. In this case, the graph $G$ has order~$n = 7$ and is obtained from a $5$-cycle by selecting two non-adjacent vertices on the cycle and adding a pendant edge to these two vertices. In this case, the set consisting of the two vertices of degree~$3$ and any two vertices of degree~$2$ is an LTD-set of $G$, implying that $\gamma_t^L(G) \le 4 < \frac{2}{3}n$, a contradiction. Hence, $G' \notin \cF_\tdom$. Thus, $\gamma_t^L(G') \le \frac{2}{3}n' = \frac{2}{3}n - 2$. Every $\gamma_t^L$-set of $G'$ can be extended to an LTD-set of $G$ by adding to it the vertices $v_2$ and $v_3$, implying that $\gamma_t^L(G) \le \frac{2}{3}n$, a contradiction.

(h) Suppose that there is a $(1,3,3,1)$-sequence in $G$. In this case, $n' = n - 4 \ge 3$ and $G'$ contains at most two components. Let $u_i$ be the neighbor of $v_i$ not on $P$ for $i \in \{2,3\}$. Possibly, $u_2 = u_3$. By parts~(a) and~(g), the vertex $u_i$ has degree~$3$ in $G$ for $i \in \{2,3\}$. Suppose that $G'$ is disconnected. In this case, by parts~(a)--(g), neither component of $G'$ belongs to $\cF_\tdom$. By linearity, we therefore have that $\gamma_t^L(G') \le \frac{2}{3}n' = \frac{2}{3}(n-4) < \frac{2}{3}n - 2$. Every $\gamma_t^L$-set of $G'$ can be extended to an LTD-set of $G$ by adding to it the vertices $v_2$ and $v_3$, implying that $\gamma_t^L(G) < \frac{2}{3}n$, a contradiction. Hence, $G'$ is connected. Recall that $n' \ge 3$. By parts~(a)--(g), we note that $G' \notin \cF_\tdom$, implying that $\gamma_t^L(G') \le \frac{2}{3}n' < \frac{2}{3}n - 2$. Every $\gamma_t^L$-set of $G'$ can be extended to an LTD-set of $G$ by adding to it the vertices $v_2$ and $v_3$, implying that $\gamma_t^L(G) < \frac{2}{3}n$, a contradiction.

Thus, the proof of the subclaim is complete.
\end{subclaimproof}

\medskip
We now return to the proof of Claim~\ref{c:claim1}. By Subclaim~\ref{c:claim1.1}(f), the neighbor of every vertex of degree~$1$ has degree~$3$ in $G$. Further by Subclaim~\ref{c:claim1.1}(a) and (g), such a vertex of degree~$3$ has both its other two neighbors of degree~$3$. Therefore the existence of a vertex of degree~$1$ implies that there is a $(1,3,3)$-sequence in $G$. In this case, $n' = n - 3 \ge 4$. Let $u_2$ be the neighbor of $v_2$ not on $P$, and let $u_3$ and $w_3$ be the two neighbors of $v_3$ not on $P$. By our earlier observations, the vertex $u_2$ has degree~$3$ in $G$, and, by Subclaim~\ref{c:claim1.1}(h),  both vertices $u_3$ and $w_3$ have degree at least~$2$ in $G$.

Suppose that $G'$ contains a component that belongs to $\cF_\tdom$. By Subclaim~\ref{c:claim1.1}, this is only possible if $u_3$ and $w_3$ are adjacent and both vertices have degree~$2$ in $G$. In this case, $G[\{v_3,u_3,w_3\}]$ is a triangle in $G$. We now consider the connected graph $G^* = G - \{v_1,v_2,v_3,u_3,w_3\}$ of order~$n^* = n - 5$. Since $u_2$ has degree~$2$ in $G^*$, we note that $n^* \ge 3$ and $G^* \notin \cF_\tdom$. Hence, $\gamma_t^L(G^*) \le \frac{2}{3}n^* = \frac{2}{3}(n-5) < \frac{2}{3}n - 3$. Every $\gamma_t^L$-set of $G^*$ can be extended to an LTD-set of $G$ by adding to it the vertices $v_2$, $v_3$, and $u_3$, implying that $\gamma_t^L(G) \le \gamma_t^L(G^*) + 3 < \frac{2}{3}n$, a contradiction. Hence, no component of $G'$ belongs to the family~$\cF_\tdom$. By linearity, we therefore have that $\gamma_t^L(G') \le \frac{2}{3}n' = \frac{2}{3}n - 2$. Every $\gamma_t^L$-set of $G'$ can be extended to an LTD-set of $G$ by adding to it the vertices $v_2$ and $v_3$, implying that $\gamma_t^L(G) < \frac{2}{3}n$, a contradiction. This completes the proof of Claim~\ref{c:claim1}.
\end{claimproof}

\medskip
By Claim~\ref{c:claim1}, every vertex in $G$ has degree~$2$ or~$3$.

\begin{claim}
\label{c:claim2}
The graph $G$ is triangle-free.
\end{claim}
\begin{claimproof} Suppose that $G$ contains a triangle $T$. Among all triangles in $G$, let $T$ contain the maximum number of vertices of degree~$2$ in $G$. Let $V(T) = \{v_1,v_2,v_3\}$, where $2 \le \deg_G(v_1) \le \deg_G(v_2) \le \deg_G(v_3) \le 3$. Since $n \ge 7$, the triangle $T$ contains at most two vertices of degree~$2$, and so $\deg_G(v_3) = 3$. Let $G' = G - V(T)$ and let $G'$ have order~$n'$, and so $n' = n - 3 \ge 4$.

Suppose that $\deg_G(v_1) = 2$. We note that $\deg_G(v_2) = 2$ or $\deg_G(v_2) = 3$. Since every vertex in $G$ has degree~$2$ or~$3$, no component of $G'$ belongs to $\cF_\tdom$. Hence by linearity, $\gamma_t^L(G') \le \frac{2}{3}n' = \frac{2}{3}n - 2$. Every $\gamma_t^L$-set of $G'$ can be extended to an LTD-set of $G$ by adding to it the vertices $v_2$ and $v_3$, implying that $\gamma_t^L(G) \leq \frac{2}{3}n$, a contradiction. Hence, $\deg_G(v_1) = 3$, implying that every vertex in $T$ has degree~$3$ in $G$. Hence by our choice of the triangle $T$, no vertex of degree~$2$ in $G$ belongs to a triangle.

Let $u_i$ be the neighbor of $v_i$ not in the triangle $T$ for $i \in [3]$. We note that the vertices $u_1$, $u_2$, and $u_3$ are not necessarily distinct. Suppose that $G'$ contains no component that belongs to $\cF_\tdom$. By linearity, this yields $\gamma_t^L(G') \le \frac{2}{3}n' = \frac{2}{3}n - 2$. Every $\gamma_t^L$-set of $G'$ can be extended to an LTD-set of $G$ by adding to it the vertices $v_2$ and $v_3$, implying that $\gamma_t^L(G) \le \frac{2}{3}n$, a contradiction. Hence, $G'$ contains a component that belongs to $\cF_\tdom$. Since $n \ge 7$ and no vertex of degree~$2$ in $G$ belongs to a triangle, this is only possible if either $G' \cong K_{1,3}$ or if $G'$ contains a $K_2$-component.

On the one hand, if $G' \cong K_{1,3}$, then the three vertices $u_1$, $u_2$, and $u_3$ are leaves in $G'$ that are adjacent to a common neighbor (of degree~$3$) in $G'$. In this case, the graph $G$ is determined and $n = 7$, and the set $V(T) \cup \{u_1\}$ is an LTD-set of $G$, implying that $\gamma_t^L(G) \le 4 < \frac{2}{3}n$, a contradiction.

On the other hand, if $G'$ contains a $K_2$-component, then renaming vertices if necessary, we may assume that $u_1$ and $u_2$ belong to such a component. We note that $u_1$ and $u_2$ both have degree~$2$ in $G$, and $u_1v_1v_2u_2u_1$ is a $4$-cycle in $G$. Further we note that in this case, $G'$ contains two components, where the second component contains the vertex $u_3$. We now consider the graph $G^* = G - \{v_1,v_2,v_3,u_1,u_2\}$. Let $G^*$ have order~$n^* = n - 5$. By the fact that $\delta(G) \geq 2$ by Claim~\ref{c:claim1}, the graph $G^* \notin \cF_\tdom$, implying that $\gamma_t^L(G^*) \le \frac{2}{3}n^* = \frac{2}{3}(n-5) < \frac{2}{3}n - 3$. Every $\gamma_t^L$-set of $G^*$ can be extended to an LTD-set of $G$ by adding to it, for example, the vertices $u_2, v_2$ and $v_3$, implying that $\gamma_t^L(G) < \frac{2}{3}n$, a contradiction.
\end{claimproof}

\medskip
By Claim~\ref{c:claim2}, the graph $G$ is triangle-free. We show next that there is no vertex of degree~$2$.

\begin{claim}
\label{c:claim3}
The graph $G$ is a cubic graph.
\end{claim}
\begin{claimproof} Suppose, to the contrary, that $\delta(G) = 2$. As before, we obtain a series of structural properties of the graph $G$ that show that certain $(d_1,d_2,\ldots,d_k)$-sequences are forbidden. These forbidden sequences will enable us to deduce the desired result of the claim that $G$ must be a cubic graph.

\begin{subclaim}
\label{c:claim3.1}
The following properties hold in the graph $G$.
\begin{enumerate}
\item[{\rm (a)}] There is no $(2,2,2)$-sequence.
\item[{\rm (b)}] There is no $(2,3,2)$-sequence.
\item[{\rm (c)}] There is no $(2,2,3)$-sequence.
\item[{\rm (d)}] There is no $(2,2)$-sequence.
\item[{\rm (e)}] There is no $(2,3,3)$-sequence.
\end{enumerate}
\end{subclaim}
\begin{subclaimproof} (a) Suppose that there is a $(2,2,2)$-sequence in $G$. In this case, $n' = n - 3 \ge 4$. Since $n \ge 7$, $\delta(G) = 2$, and $G$ contains no triangle, no component of $G'$ belongs to $\cF_\tdom$. Hence by linearity, $\gamma_t^L(G') \le \frac{2}{3}n' = \frac{2}{3}n - 2$. Every $\gamma_t^L$-set of $G'$ can be extended to an LTD-set of $G$ by adding to it the vertices $v_2$ and $v_3$, implying that $\gamma_t^L(G) \le \frac{2}{3}n$, a contradiction.

(b) Suppose that there is a $(2,3,2)$-sequence in $G$. As before, $n' = n - 3 \ge 4$. Suppose that $G'$ contains a component that belongs to $\cF_\tdom$. Since there is no $(2,2,2)$-sequence and $n \ge 7$, and since $\delta(G) \ge 2$ and $G$ contains no triangle, this is only possible if $G' \cong K_{1,3}$. But then the graph $G$ is determined and $n = 7$, and $\gamma_t^L(G') = 4 < \frac{2}{3}n$ (by considering the set $N[v_2]$), a contradiction. Hence, no component of $G'$ belongs to $\cF_\tdom$. By linearity, this yields $\gamma_t^L(G') \le \frac{2}{3}n' = \frac{2}{3}n - 2$. Every $\gamma_t^L$-set of $G'$ can be extended to an LTD-set of $G$ by adding to it the vertices $v_2$ and $v_3$, implying that $\gamma_t^L(G) \le \gamma_t^L(G') + 2 \le \frac{2}{3}n$, a contradiction.

(c) Suppose that there is a $(2,2,3)$-sequence in $G$. Since there is no $(2,2,2)$-sequence and no $(2,3,2)$-sequence in $G$, every vertex different from $v_2$ that is adjacent to $v_1$ or $v_3$ has degree~$3$ in $G$. Together with our earlier observations, the graph $G'$ therefore cannot contain a component that belongs to $\cF_\tdom$. By linearity, we have $\gamma_t^L(G') \le \frac{2}{3}n' = \frac{2}{3}n - 2$. As before this yields $\gamma_t^L(G) \le \gamma_t^L(G') + 2 \le \frac{2}{3}n$, a contradiction.

(d) Since there is no $(2,2,2)$-sequence and no $(2,2,3)$-sequence, there can be no $(2,2)$-sequence in $G$ noting that every vertex has degree~$2$ or~$3$.

(e)  Suppose that there is a $(2,3,3)$-sequence in $G$. Since there is no $(2,2,2)$-sequence, no $(2,3,2)$-sequence, and no $(2,2,3)$-sequence in $G$, the graph $G'$ cannot contain a component that belongs to $\cF_\tdom$. By linearity, this yields $\gamma_t^L(G') \le \frac{2}{3}n' = \frac{2}{3}n - 2$. Every $\gamma_t^L$-set of $G'$ can be extended to an LTD-set of $G$ by adding to it the vertices $v_2$ and $v_3$, implying that $\gamma_t^L(G) \le \gamma_t^L(G') + 2 \le \frac{2}{3}n$, a contradiction.

Thus, the proof of the subclaim is complete.
\end{subclaimproof}

\medskip
By Subclaim~\ref{c:claim3.1}(d), there is no $(2,2)$-sequence. Hence every vertex of degree~$2$ has both its neighbors of degree~$3$. Moreover since there is no $(2,3,2)$-sequence, every vertex of degree~$3$ has at most one neighbor of degree~$2$. But this would imply the existence of a $(2,3,3)$-sequence, contradicting Subclaim~\ref{c:claim3.1}(e). Therefore, there can be no vertex of degree~$2$ in $G$, that is, $G$ is a cubic graph. This completes the proof of  Claim~\ref{c:claim3}.
\end{claimproof}

\medskip
By Claim~\ref{c:claim3}, the graph $G$ is a cubic graph. Recall that $G$ is triangle-free. We now consider a $(3,3,3)$-sequence. The graph $G'$ cannot contain a component that belongs to $\cF_\tdom$. By linearity, this yields $\gamma_t^L(G') \le \frac{2}{3}n' = \frac{2}{3}n - 2$. Every $\gamma_t^L$-set of $G'$ can be extended to an LTD-set of $G$ by adding to it the vertices $v_2$ and $v_3$, implying that $\gamma_t^L(G) \le \gamma_t^L(G') + 2 \le \frac{2}{3}n$, a contradiction. This completes the proof of Theorem~\ref{t:thm-subcubic}.
\end{proof}

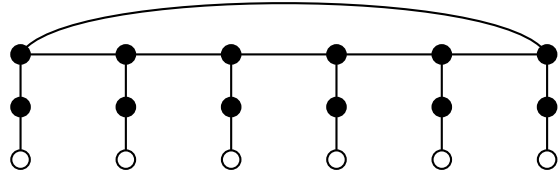
\begin{figure}[htb]
\begin{center}
\begin{tikzpicture}[scale=.7,style=thick,x=1cm,y=1cm]
\def\vr{5pt}
\path (0,0) coordinate (v1);
\path (0,1) coordinate (v2);
\path (0,2) coordinate (v3);
\draw (v1) -- (v2);
\draw (v2) -- (v3);
\path (2,0) coordinate (u1);
\path (2,1) coordinate (u2);
\path (2,2) coordinate (u3);
\draw (u1) -- (u2);
\draw (u2) -- (u3);
\path (4,0) coordinate (w1);
\path (4,1) coordinate (w2);
\path (4,2) coordinate (w3);
\draw (w1) -- (w2);
\draw (w2) -- (w3);
\path (6,0) coordinate (x1);
\path (6,1) coordinate (x2);
\path (6,2) coordinate (x3);
\draw (x1) -- (x2);
\draw (x2) -- (x3);
\path (8,0) coordinate (y1);
\path (8,1) coordinate (y2);
\path (8,2) coordinate (y3);
\draw (y1) -- (y2);
\draw (y2) -- (y3);
\path (10,0) coordinate (z1);
\path (10,1) coordinate (z2);
\path (10,2) coordinate (z3);
\draw (z1) -- (z2);
\draw (z2) -- (z3);
\draw (v3) -- (u3);
\draw (u3) -- (w3);
\draw (w3) -- (x3);
\draw (x3) -- (y3);
\draw (y3) -- (z3);
\draw (v1) [fill=white] circle (\vr);
\draw (v2) [fill=black] circle (\vr);
\draw (v3) [fill=black] circle (\vr);
\draw (u1) [fill=white] circle (\vr);
\draw (u2) [fill=black] circle (\vr);
\draw (u3) [fill=black] circle (\vr);
\draw (w1) [fill=white] circle (\vr);
\draw (w2) [fill=black] circle (\vr);
\draw (w3) [fill=black] circle (\vr);
\draw (x1) [fill=white] circle (\vr);
\draw (x2) [fill=black] circle (\vr);
\draw (x3) [fill=black] circle (\vr);
\draw (y1) [fill=white] circle (\vr);
\draw (y2) [fill=black] circle (\vr);
\draw (y3) [fill=black] circle (\vr);
\draw (z1) [fill=white] circle (\vr);
\draw (z2) [fill=black] circle (\vr);
\draw (z3) [fill=black] circle (\vr);
\draw (v3) to[out=60,in=120, distance=1.5cm ] (z3);
\end{tikzpicture}
\end{center}
\begin{center}
\caption{The $2$-corona $C_6 \circ P_2$ of a $6$-cycle.}
\label{2corona-cycle}
\end{center}
\end{figure}

For $k \ge 3$, the $2$-corona $G = C_k \circ P_2$ of a cycle $C_k$ has order~$n = 3k$ and by Observation~\ref{obs:conj-tight}, it has locating total domination number equal to its total domination number, that is, $\gamma_t^L(G) = \gamma_t(G) = 2k = \frac{2}{3}n$. See Figure~\ref{2corona-cycle} for an illustration. Moreover for $k \ge 1$, the $2$-corona $G = P_k \circ P_2$ of a path $P_k$ has order~$n = 3k$ and also satisfies $\gamma_t^L(G) = \gamma_t(G) = 2k = \frac{2}{3}n$. Thus, we obtain the following.

\begin{proposition}
\label{tightness-of-subcubic}
There are infinitely many connected twin-free subcubic graphs $G$ of order $n$ with $\gamma_t^L(G) = \frac{2}{3}n$.
\end{proposition}

\section{Conclusion}\label{sec:conclu}

We have proved Conjecture~\ref{conj} for several important graph classes: cobipartite graphs, split graphs, block graphs and subcubic graphs.

It would be interesting to extend these results to larger classes, for example chordal graphs (which include split graphs and block graphs). Another interesting subclass of chordal graphs to consider is the class of interval graphs.

It would also be interesting to prove that the bound $\gLT(G)\leq\frac{n}{2}$ holds for sufficiently large (twin-free) connected cubic graphs, as conjectured in~\cite{hl12}.

\end{document}